\documentclass[11pt,a4paper,twoside,openright]{article}
\usepackage[english]{babel}
\usepackage{amsmath,amssymb,amsfonts,epic}
\usepackage[dvips]{graphicx}
\usepackage[latin1]{inputenc} 
\usepackage{chicago}

\usepackage{setspace}
\newcommand{\e}{\varepsilon}
\newcommand{\N}{\mathbb{N}}
\newcommand{\R}{\mathbb{R}}
\newcommand{\A}{\mathcal{A}}                         % alphabet

\newcommand{\pr}{\mathbb{P}_\theta}             % proba sous theta
\newcommand{\pru}{\mathbb{P}_{\theta_1}}             % proba sous theta1
\newcommand{\prd}{\mathbb{P}_{\theta_2}}             % proba sous theta2
\newcommand{\esp}{\mathbb{E}_\theta}          % esperance sous theta
\newcommand{\pro}{\mathbb{P}_{0}}                % proba sous theta0
\newcommand{\espo}{\mathbb{E}_{0}}             % esperance sous theta0
\newcommand{\Q}{Q_\theta}
\newcommand{\Qb}{Q_{\theta (\beta)}}                            % Q sous theta
\newcommand{\Qu}{Q_{\theta_1}}                               % Q sous theta
\newcommand{\Qd}{Q_{\theta_2}}                               % Q sous theta
\newcommand{\Qo}{Q_{\theta_0}}                      % Q sous theta0

\newcommand{\1}[1]{1\! \mathrm{l}\{ #1\}}         % indicateur d'un ensemble
\newcommand{\E}{\mathcal{E}}                       % espace des trajectoires

\newtheorem{thm}{Theorem}
\newtheorem{defi}{Definition}

\newtheorem{cor}{Corollary}
\newtheorem{assumption}{Assumption}
\newtheorem{lemma}{Lemma}

\newenvironment{proof}[1][Proof]{\textbf{#1} }{\ \rule{0.5em}{0.5em}}

\title{Parameter estimation in pair hidden Markov models.} 
\author{Ana Arribas-Gil $^\star $, Elisabeth Gassiat $^\star $, Catherine
Matias$^{\star\star} $}

\begin{document}

\thispagestyle{empty}
\maketitle

\begin{center}
$^\star$ UMR CNRS 8628.  \'Equipe Probabilit\'es, Statistique et
 Mod\'elisation, Universit\'e Paris-Sud, 
B\^atiment 425, Universit\'e de Paris-Sud, 91405 Orsay Cedex,
 France. \\ E-mail: \{ana.arribasgil, elisabeth.gassiat\}@math.u-psud.fr \\
 $^{\star\star}$ UMR CNRS 8071. Laboratoire Statistique et G\'enome, 
  Tour Evry  2, 523  pl.   des Terrasses  de l'Agora,  91~000 Evry,  France.
E-mail:  matias@genopole.cnrs.fr
\end{center}

\begin{abstract}
This paper deals with parameter estimation in pair hidden Markov models 
(pair-HMMs). We first 
provide a rigorous formalism for these models and discuss possible 
definitions of likelihoods. The model being biologically motivated, some 
restrictions with respect to the full parameter space naturally occur. 
Existence of two 
different Information divergence rates is established and divergence 
property (namely positivity at values different from the true one) is 
shown under additional assumptions. This yields consistency for the parameter 
in parametrization schemes for which the divergence property 
holds. Simulations illustrate different cases which are not covered by our 
results. 
\end{abstract}

%\noindent {\it AMS Classification: 62F12, 62M09, 62B10.}
%{\small }

%% AMS Classif
%% 62F : parametric inference / 62F12 : asymptotic properties
%% 62M : inference from  stochastic processes/ 62M09 : Non-Markovian processes
%% estimation
%% 62B10 : information theoretic topics

\section{Introduction}
\subsection{Background}

\begin{figure}[!h]
\begin{center}
\setlength{\unitlength}{1mm}
\begin{picture}(60,60)
%Axes
\put(6,6){\line(1,0){46}} \put(6,6){\line(0,1){46}}
\multiput(16,6)(10,0){4}{\dashbox{1}(0,46)}
\multiput(6,16)(0,10){4}{\dashbox{1}(46,0)}
%Lettres
\put(9,0){\Large A}\put(19,0){\Large A}\put(29,0){\Large
T}\put(39,0){\Large G} \put(0,9){\Large C}\put(0,19){\Large
T}\put(0,29){\Large G}\put(0,39){\Large G}
%Chemin
\thicklines
\drawline[5](6,6)(16,16)(26,16)(36,26)(46,36)(46,46)
\linethickness{0.35mm}
\put(16,16){\line(1,0){10}}\put(46,36){\line(0,1){10}}
\end{picture}
\end{center}
\caption{Graphical  representation  of  an  alignment  between  two  sequences
  $X=AATG$ and $Y=CTGG$. The displayed alignment is $ 
  \stackrel{A}  {\text{\it  \scriptsize C}}  \stackrel{A}{\text{  \large -  }}
  \stackrel{T}{\text{\it  \scriptsize T}} \stackrel{G}{  \text{\it \scriptsize
      G}}   \stackrel{\text   {  \large   -   }}{\text{\it  \scriptsize   G}}.
  $}\label{grid} 
\end{figure}
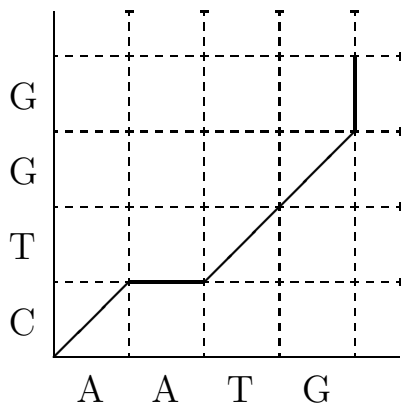

Sequence   alignment  has   become  one   of  the   most  powerful   tools  in
bioinformatics. Biological sequences are  aligned for instance (and among many
other  examples)  to infer  gene  functions,  to  construct or  use  protein
databases  or to  construct phylogenetic  trees. Concerning  this  last topic,
current methods first  align the sequences and then  infer the phylogeny given
this  fixed alignment.  This  approach contains  a  major flaw  since the  two
problems are largely intertwined.
Indeed, the alignment problem consists in
retrieving the places, in the observed sequences, where
substitution/deletion/insertion  events  have occurred,  due  to the  evolution
process.
In the pair alignment problem,
the  observations  consist in  a  couple  of sequences  $X_{1:n}=X_1\ldots  X_n$  and
$Y_{1:m}=Y_1 \ldots  Y_m$ with values on  a finite state  alphabet $\A$
($\A=\{A,C,G,T\}$ for DNA sequences). It is assumed that
the sequences share a common ancestor. According to biological evolution,
the sequence of the ancestor evolves and letters in each site may
change (substitution event), or be deleted (deletion event), or new letters
may be inserted in the  sequence (insertion event). This process finally leads
to the two different observed sequences.
%When facing to sequences having evolved from a common ancestor, the alignment problem consists in retrieving
% the places where substitution/deletion/insertion events have occurred.
A most convenient  way of displaying alignments is  a graphical
representation as a  path through a rectangular grid (see Figure
\ref{grid}). A diagonal move corresponds to a match  between the two
sequences, whereas horizontal and vertical moves  correspond to
insertion-deletion events. This path consists of
steps  $\e_t$, $t=1,\ldots,l$, where $\e_t$
represents either a match  ($\e_t =(1,1)$) or an insertion-deletion
event ($\e_t =(1,0)$ or $(0,1)$). The length of the alignment is $l$, and
satisfies
 \begin{equation} \label{longueur}
n\vee m \leq l \leq n+m .
\end{equation}
Here $n\vee m$ denotes the maximum value between $n$ and $m$.
The multiple alignment problem is the same, except that one has to retrieve
the places where substitution/deletion/insertion events have occurred on the basis
of a set of (more than two) sequences. 

Aligning  two sequences relies on  the choice of  a score optimization
scheme  (for instance,  the Needleman-Wunsch  algorithm  \cite{Needleman}) and
therefore  the obtained alignments  depend on  the score  parameters. Choosing
these  score  parameters  in the  most  objective  way  appears as  a  crucial
issue.  Because  evolution  is  the  force that  promotes  divergence  between
biological sequences, it is desirable  to consider biological alignment in the
context of evolution.
Now, given an evolution model,  optimal choices of the score parameters depend
on  the underlying  unknown mutation  rates and  thus on  the phylogeny  to be
inferred after the alignment. The  existence of such a vicious circle explains
the  emergence of  probabilistic models  where optimal  alignment  and evolution  parameters
estimation are achieved at the same time.\\

Relying on  a pioneering  work by Bishop  and Thompson  \citeyear{Bishop}, Thorne,
Kishino  and  Felsenstein \citeyear{TKF}  were  the  first  to provide  a  maximum
likelihood approach to  the alignment of a pair of DNA  sequences based on a
rigorous model of sequence evolution (referred to as the TKF model).
This model  has become quite classical  nowadays. In this setup,  each site is
independently hit  by a substitution  or deleted, and insertions  occur between
two sites or at both ends of the sequence. Each one of those events occurs at a
specific rate.  When  a substitution or an insertion  occurs, a new nucleotide
is  drawn randomly  according to  some probability  distribution on  the state
space $\{A,C,G,T\}$.
One of the advantages of the TKF model lies in its exact correspondence with a
model containing a hidden Markov structure, ensuring the existence of powerful
algorithmic tools  based on dynamic  programming methods. More  precisely, the
TKF  evolution model  fits into  the  concept of  a pair  hidden Markov  model
(pair-HMM), as first formally described in \cite{Durbin}. \\

Observations  in a  pair-HMM   are formed  by a
couple  of sequences  (the  ones to  be  aligned) and  the  model assumes  that
the  hidden (i.e. non observed) alignment sequence $\{\e_t \}_{t}$ is a Markov
chain that determines the probability distribution of the observations.
Since the seminal paper \cite{TKF}, an abundant literature aroused in which 
parameter  estimation  occurs  in   a  pair-HMM.  Thorne,  Kishino  and
Felsenstein \citeyear{TKF2}  slightly improved their  original model to  take into
account  insertion and  deletion  of  entire fragments  (and  not only  single
nucleotides).   The  TKF  model  approaches  have been  further  developed  in
\cite{HeinWiuf,Metzler,KnudsenMiyamoto,Miklos}, for
instance.  Let us  also mention  that  pair-HMMs were  recently combined  with
classical hidden Markov models (HMMs) for {\it ab
initio} prediction of genes \cite{MeyerDurbin,Pachter,Hobolth}.\\

The main difference between pair-HMMs and classical HMMs lies in the
observation of a {\it pair} of sequences instead of a {\it single}
one.  From a  practical  point of  view, the  two  above models  are
not  very different  and  classical  algorithms  such  as 
forward  or    Viterbi algorithms are still valid and efficient
in the pair-HMM context (we refer to \cite{Durbin} for a complete
description of those techniques). Forward  algorithm  allows to
compute  the  likelihood  of the  two  observed sequences and thus,
by means of a maximisation  technique, to approximate the maximum
likelihood  estimator (MLE) of the  parameters. Numerical
maximisation approaches are  commonly used (\cite{TKF}) but
statistical  approaches using the Expectation-Maximisation  (EM)
algorithm  and  its variants  (Stochastic  EM, Stochastic
Approximation EM) have  recently been  explored \cite{Holmes,Ana}.
Viterbi algorithm is designed to  reconstruct the most probable
hidden path, thus giving the alignment. From a Bayesian
point of view, it is also interesting to provide a posterior
distribution for parameters and alignments. This can be done
with MCMC procedures  needing again the use of  forward
algorithm \cite{Metzler,Ana}. \\

Nonetheless, from a theoretical point
of view, pair-HMMs and classical HMMs are completely different. In
particular, up to our knowledge,
there is no theoretical proofs that the maximum likelihood
procedure  nor the Bayesian estimation  give consistent
estimators of the pair-HMM parameters (though it is the case for instance
for regular HMMs with finite state space, see \cite{BaumPetrie} concerning MLE
consistency;  see also  \cite{Caliebe} for  the convergence  of the  maximum a
posteriori hidden path).

This paper  is thus concerned with  statistical properties of  parameter estimation
procedures in  pair-HMMs. \\

\subsection{Roadmap}
In Section 2, the pair-HMM  is described, together with some properties
of the distribution of observed sequences. Then we state possible likelihood functions,
to be compared with the criterion that
is optimized in pair-HMM algorithms. We then interpret this last one  as a likelihood function.
\\
To investigate consistency of estimators obtained by maximization, one has
to understand the asymptotic behaviour of the criteria.  We adopt the Information Theory
terminology and call {\it Information divergence rates} the difference between
the limiting  values of  the log-likelihoods at  the (unknown)  true parameter
value and at another parameter value.
Indeed, the general model described below may be interpreted as a channel
transmitting the input $X_{1:n}$ with possible errors, insertions or deletions, leading
to  the output  $Y_{1:m}$ (see  for instance  \cite{Davey,Levenshtein}  on the
topic of error correcting codes and also
\cite{Csiszar,Cover} for a general introduction to Information Theory).
In this setting, {\it Information divergence rates} have a precise meaning (in terms of
coding or transmission qualities). In a statistical setting such as ours, they are
interpreted as  divergences that should have a unique minimum at the true parameter
value (divergence property).
Section 3 is devoted
to  the  existence  and  properties  of such  limit  functions  (see  Theorems
\ref{thdivergence} and \ref{contrast}).\\
Section 4,  then, gives  the statistical consequences  in terms  of consistent
estimation of  the parameters  obtained via MLE  or Bayesian  estimation using
pair-HMM algorithms (see Theorems \ref{thcons}, \ref{thbayes}).
According to these results, consistency holds for the parameter in parametrization schemes
for which the divergence property holds for the associated Information divergence rate. \\
In a last section, we present several simulation results to investigate situations
in which the divergence property is not established.
We illustrate the consistency results in cases where Theorem
\ref{thcons} applies, as may be seen on
numerical computations of information divergence rates.
We also compare the limiting values
of different criteria and give some interpretations. Unfortunately, despite the
positive results that we obtain we are not yet in terms of completely
validate pair-HMM algorithms in every situation. 

\section{The pair hidden Markov model}
\subsection{Model description}
We now describe in details the pair-HMM. Consider a
stationary ergodic Markov  chain $\{\e_t\}_{t \geq 1}$ on the state
space $\E= \{ (1,0); (0,1) ; (1,1)\}$, with transition matrix $\pi$ and
stationary distribution $\mu=(p,q,r)$. This chain generates a random
walk $\{Z_t\}_{t \geq 0}$ with values in the two-dimensional integer
lattice $\N \times \N$, by letting $Z_0 =(0,0)$ and $Z_t  =
\sum_{1\leq s\leq t} \e_s$. The  coordinate random variables
corresponding to  $Z_t$ at time $t$ are denoted  by  $(N_t,M_t)$
({\it i.e.} $Z_t=(N_t,M_t)$). We  shall   either  use  the  notation
$\pi(\e_s,\e_{s+1})$   to  denote  the transitions  probabilities of
the  matrix  $\pi$,  or  explicit symbols  like $\pi_{HV}$
indicating a transition from  state $H=(1,0)$  to  state  $V=(0,1)$
($H$ stands  for  {\it
  horizontal} move, $V$ for {\it vertical} move and $D=(1,1)$ for {\it diagonal}  move).

Conditional on  the hidden random walk, the  observations are drawn
according to the following scheme. At  time $t$, if $\e_t=(1,0)$ then a random
variable $X$ is drawn (emitted) according  to some probability distribution $f$ on $\A$,
if  $\e_t=(0,1)$  then  a random  variable  $Y$  is  drawn (emitted) according  to  some
probability  distribution $g$  on $\A$  and  finally, if  $\e_t=(1,1)$ then  a
couple  of random  variables $(X,Y)$  is drawn (emitted) according to  some probability
distribution $h$ on $\A\times \A$. Conditionally to the hidden Markov chain
$\{\e_t\}_{t \geq 1}$, all emitted random variables are independent.
This  model is  described by  the  parameter $\theta=(\pi,  f,g,h) \in  \Theta
$.
The conditional distribution of the observations thus writes
\begin{multline} \label{conditional}
\pr(X_{1:N_t}, Y_{1:M_t} | \e_{1:t}, \{\e_s\}_{s>t}, \{X_i, Y_j\}_{i\neq N_s, j\neq M_s, 0\leq s\leq t})=
\pr(X_{1:N_t}, Y_{1:M_t} |  \e_{1:t} ) \\
= \prod_{s=1}^t f(X_{N_s})^{\1{\e_s=(1,0)}}
g(Y_{M_s})^{\1{\e_s=(0,1)}} h(X_{N_s},Y_{M_s})^{\1{\e_s=(1,1)}} ,
\end{multline}
where $\1{\cdot}$  stands for the  indicator function.
Moreover,  the complete distribution $\pr$ is given by
$$
\pr(\e_{1:t},  X_{1:N_t},   Y_{1:M_t})  =  \mu(\e_1)\Big   \{  \prod_{s=2}^{t}
\pi(\e_{s-1}, \e_s) \Big\} \pr(X_{1:N_t}, Y_{1:M_t} | \e_{1:t}) .
$$
Here we denote by $\pr$  (and $\esp$) the  induced probability distribution
(and corresponding  expectation) on $\E^{\N}  \times \A^{\N} \times
\A^{\N}$ and  $\theta_0$ the true parameter corresponding  to the distribution
of the observations (we shall abbreviate to $\pro$ and $\espo$ the probability
distribution and expectation under parameter $\theta_0$).
Note that a necessary condition for identifiability of the parameter $\theta$ is that the occurrence probability of two aligned letters differs from the product probabilities of these letters. That is:
\begin{assumption}\label{assum:ident}
 $$\exists x, y \in \A , \text{ such that } h(x,y) \neq f(x) g(y).$$
\end{assumption}
Indeed, if  $h=f g$, then \eqref{conditional} gives
$$
\pr(X_{1:N_t},  Y_{1:M_t}  | \e_{1:t}  )  =\left\{ \prod_{i=1}^{N_t}  f(X_{i})
\right\}   \left\{\prod_{j=1}^{M_t}  g(Y_{j})   \right  \}   =  \pr(X_{1:N_t},
Y_{1:M_t} ) .
$$
Thus,  in  this  case,  the  distribution  of the  observations  is
independent  from  the  hidden  process  and the  parameter  $\pi$  cannot  be
identified.  In the following, we shall always work under Assumption~\ref{assum:ident}.

\subsection{Observations and likelihoods.}
Statisticians define log-likelihoods to be functions of the parameter,
that are equal to the logarithm of the probability of the
observations. Here, to state what log-likelihoods are, one has to decide
what do the observed sequences $(X_{1:n}, Y_{1:m})$ represent.  Indeed, one may
interpret it in at least two different ways:
\begin{itemize}
\item[(a)]
It is the  observation of emitted sequences until some time  $t$, so that the
log-likelihood should be $\log \pr(X_{1:N_t}, Y_{1:M_t})$.
Here, the probability is that
of the observed sequences {\it and} a point of the hidden process $Z_t=(N_t,M_t)$;
\item[(b)]
Each observed sequence is one of the emitted sequences $X_{1:N_t}$ for some $t$
and $Y_{1:M_s}$ for some $s$, knowing nothing on the hidden process (that is whether $t=s$,
or $t>s$, or $t<s$),
so that the log-likelihood should be $\log \pr(X_{1:n}, Y_{1:m})$. Here, the probability
is the marginal distribution of the sequences.
\end{itemize}
It should be  now noted that none  of those quantities is the  one computed by
pair-HMM algorithms.  We will come back  to this fact later. Note also that we
imposed the true underlying alignment to pass through the fixed point $(0,0)$
(namely, we assumed $Z_0=(0,0)$) which is  not the more general setup (and may
introduce  a  bias  in  practical  applications). However,  we  restrict  our
attention to this particular setup.

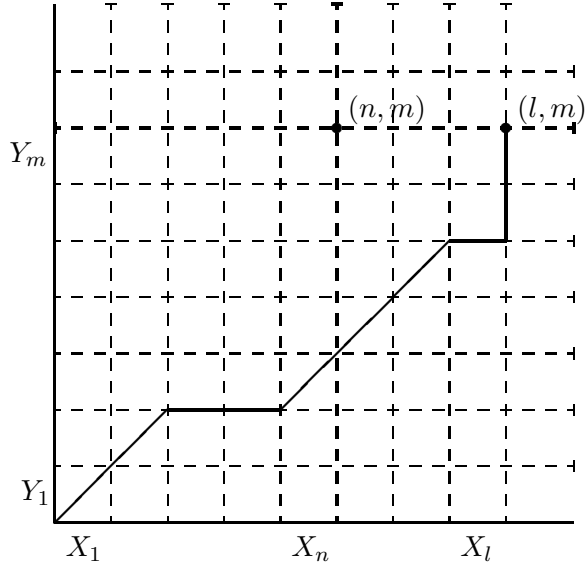
\begin{figure}[!h]
\begin{center}
\setlength{\unitlength}{1.5mm}
\begin{picture}(60,60)
%Axes
\put(6,6){\line(1,0){46}} \put(6,6){\line(0,1){46}}
\multiput(11,6)(5,0){8}{\dashbox{1}(0,46)}
\multiput(6,11)(0,5){8}{\dashbox{1}(46,0)}
%Lettres
\put(7,3){$X_1$}\put(27,3){$X_n$}\put(42,3){$X_l$}
\put(3,8){$Y_1$}\put(2,38){$Y_m$} \put(31,41){\circle*{1}}
\put(32,42){$(n,m)$} \put(46,41){\circle*{1}} \put(47,42){$(l,m)$}
%Chemin
\thicklines
\drawline[5](6,6)(16,16)(26,16)(36,26)(41,31)(46,31)(46,36)
\linethickness{0.35mm} \put(31,6){\dashbox{1}(0,46)}
\put(6,41){\dashbox{1}(46,0)}
\put(16,16){\line(1,0){10}}\put(46,36){\line(0,1){5}}
\put(41,31){\line(1,0){5}}\put(46,31){\line(0,1){5}}
\end{picture}
\end{center}
\caption{Graphical representation of an alignment of sequences
$X_{1:n}$ and $Y_{1:m}$ not passing through the point
$(n,m)$.}\label{grid2}
\end{figure}

First, we introduce some notations to make the previous quantities more precise.\\
Let us  consider the set $\E_\infty$ of  all the possible trajectories  of the hidden
path and the set $\E_{n,m}$ of trajectories passing through the point $(n,m)$:
 \begin{eqnarray}
&  \E_\infty   &=  \{(0,1);  (1,0);  (1,1)\}^\N   =  \{e=(e_1,  e_2,   \ldots)  \}  =
\E^{\N} \label{traject} , \\
  &\E_{n,m} &= \{ e \in  \{(0,1); (1,0); (1,1)\}^l ; n\vee m \leq l \leq n+m ; \sum_{i=1}^l e_i =(n,m) \} \label{trajectnm} .
 \end{eqnarray}
The length  of  any trajectory  $e  \in \E_{n,m}$  is denoted  by $|e|$. Then, we have the following equations
\begin{equation}
\label{Pnm}
\pr(X_{1:n}, Y_{1:m}) = \sum_{e \in \E_\infty} \pr(\e_{1:\infty} =e_{1:\infty}, X_{1:n}, Y_{1:m}),
\end{equation}
\begin{equation}
\label{Pt}
\pr(X_{1:N_t}, Y_{1:M_t}) = \pr(X_{1:N_t}, Y_{1:M_t}, Z_t) = \sum_{e \in \E_{N_t,M_t} ; |e| =t} \pr(\e_{1:t} =e_{1:t}, X_{1:N_t}, Y_{1:M_t})
.
\end{equation}
As Equation \eqref{Pnm} shows, if one uses the marginal distributions as likelihood,
it means that  when observing two
sequences $X_{1:n}$ and $Y_{1:m}$, it is not assumed
that the hidden process passes through the observed point $(n,m)$.
This results in an alignment with not necessarily bounded length
(see Figure \ref{grid2}).

We shall now detail Equation \eqref{Pnm} according to possible alignments.
Among all  the trajectories in $\E_\infty$,  we shall distinguish  the ones in
$\E_{n,m}$ and the ones belonging to some set $\E_{n,p}$ (with $p>m$) or
$\E_{p,m}$ (with  $p>n$). Those last ones  need to be constrained  in order to
avoid multiple counting.
Let us  denote by $\E_{n,m}^{-H}$  (resp. $\E_{n,m}^{-V}$) the  restriction of
the   set  $\E_{n,m}$   to  trajectories   not  ending   with   an  horizontal
(resp. vertical) part. More precisely,
\begin{eqnarray*}
  \E_{n,m}^{-H}&=&\{  e=(e_1, \ldots, e_{|e|})  \in \E_{n,m}  ; \text{  If for
    some   } s, \\
&&  \hspace{4cm} e_{|e|}=e_{|e|-1}=  \ldots =e_{|e|-s+1}=(1,0)  \text{  then }
s=0\} ,\\ \E_{n,m}^{-V}&=&\{ e=(e_1, \ldots, e_{|e|}) \in \E_{n,m} ; \text{ If
  for some } s, \\
&&  \hspace{4cm} e_{|e|}=e_{|e|-1}=  \ldots =e_{|e|-s+1}=(0,1)  \text{  then }
s=0\} .
\end{eqnarray*}
These  notations  allow to  express  the  marginal distribution  $\pr(X_{1:n},
Y_{1:m})$  as a  sum  over  three different path types.
\begin{multline*}
\pr(X_{1:n}, Y_{1:m} ) =\sum_{e \in \E_{n,m}} \pr(\e_{1:|e|}=e, X_{1:n},
Y_{1:m}) \\
+ \sum_{p>n} \sum_{e\in \E_{p,m}^{-H}} \sum_{x_{n+1:p}} \pr(\e_{1:|e|}=e, X_{1:n},
X_{n+1:p}=x_{n+1:p} , Y_{1:m}) \\
+ \sum_{p>m} \sum_{e\in \E_{n,p}^{-V}} \sum_{y_{m+1:p}} \pr(\e_{1:|e|}=e, X_{1:n},
 Y_{1:m}, Y_{m+1:p}=y_{m+1:p}) .
\end{multline*}
This  form
%enlightens the  difference  between  the quantities  $\pr(X_{1:n},
%Y_{1:m})$ and $\Q(X_{1:n}, Y_{1:m})$
may not be used for the computation of the
marginal distribution $\pr(X_{1:n}, Y_{1:m})$.

We  now   give  some  recursion   formulas  that  could  lead   to  practical
implementations of this last quantity.
For any state $e \in \E$, define $\pr^{e}$ as the
distribution induced by $\pr$ conditional on $\e_1=e$. Let us also denote by
$h_X$ (resp. $h_Y$)  the marginal with respect to the first  (resp. second) coordinate of
the distribution $h$.

\begin{lemma} \label{lem:prec}
For any $n \geq 1, m\geq1,$
\begin{equation}
\label{prec}
\pr(X_{1:n},Y_{1:m} ) =  p \,  \pr^{(1,0)}(X_{1:n},Y_{1:m}  )  +q \,
\pr^{(0,1)}    (X_{1:n},Y_{1:m})     +    r   \,   \pr^{(1,1)}(X_{1:n},Y_{1:m} ) ,
\end{equation}
with the following recursions
\begin{eqnarray*}
  \pr^{(1,0)}(X_{1:n},Y_{1:m}       )&=&        f(X_1)       \{       \pi_{HH}
  \pr^{(1,0)}(X_{2:n},Y_{1:m} ) +\pi_{HV} \pr^{(0,1)}(X_{2:n},Y_{1:m} ) \\
&&  +\pi_{HD} \pr^{(1,1)} (X_{2:n}, Y_{1:m}) \} \\
\pr^{(0,1)} (X_{1:n},Y_{1:m}) &=&   g(Y_1)       \{       \pi_{VH}
  \pr^{(1,0)}(X_{1:n},Y_{2:m} ) +\pi_{VV} \pr^{(0,1)}(X_{1:n},Y_{2:m} ) \\
&& +\pi_{VD} \pr^{(1,1)} (X_{1:n}, Y_{2:m}) \}\\
 \pr^{(1,1)}(X_{1:n},Y_{1:m} ) &=&   h(X_1,Y_1)       \{       \pi_{DH}
  \pr^{(1,0)}(X_{2:n},Y_{2:m} ) +\pi_{DV} \pr^{(0,1)}(X_{2:n},Y_{2:m} ) \\
&& +\pi_{DD} \pr^{(1,1)} (X_{2:n}, Y_{2:m}) \}
\end{eqnarray*}
and initializations:
\begin{eqnarray*}
\pr^{(1,0)}(X_{1})  &=&  f(X_1) ,\quad  \pr^{(0,1)}(Y_{1})  =  g(Y_1) ,  \quad
\pr^{(1,1)} (X_1,Y_1) = h(X_1,Y_1) ,\\
\pr^{(0,1)} (X_{1:n}) &=& \frac{1}{1-\pi_{VV}}  \{ \pi_{VH}\; f(X_1) \pr^{(1,0)}
 (X_{2:n}) +\pi_{VD} \; h_X(X_1) \pr^{(1,1)} (X_{2:n})\},\\
 \pr^{(1,0)}(Y_{1:m}) &=& \frac{1}{1-\pi_{HH}} \{ \pi_{HV} \; g(Y_1) \pr^{(0,1)}
 (Y_{2:m}) +\pi_{HD} \; h_Y(Y_1) \pr^{(1,1)} (Y_{2:m})\}.
\end{eqnarray*}
\end{lemma}
Proof of Lemma~\ref{lem:prec} is trivial and therefore omitted.\\

Interpretation (a)  leads to define
the log-likelihood $\ell_t(\theta)$ as
\begin{equation}
\label{lt}
\ell_t  (\theta)  =  \log  \pr  (X_{1:N_t},  Y_{1:M_t}) ,  \quad  t  \geq  1 .
\end{equation}
But since  the underlying process $\{Z_t\}_{t \geq  0}$ is not observed,
the   quantity  $\ell_t(\theta)$  is   not  a   measurable  function   of  the
observations. More precisely, the {\it  time $t$} at which observation is made
is not observed itself. Though, if one decides
to use interpretation (a),  namely that $(X_{1:n},Y_{1:m})$ corresponds to the
observation  of  the  emitted sequences  at  a  point  of the  hidden  process
$Z_t=(N_t,M_t)$ and some {\it unknown} time $t$, one does not
use $\ell_t(\theta)$ as a log-likelihood, but rather
\begin{equation}
\label{qt}
w_t (\theta) =  \log \Q (X_{1:N_t}, Y_{1:M_t}) , \quad t \geq 1
\end{equation}
where for any integers $n$ and $m$
\begin{equation}
\Q(X_{1:n}, Y_{1:m})
=\pr ( \exists s \geq 1, Z_s =(n,m) ; X_{1:n}, Y_{1:m}).
\end{equation}
In other words, $\Q$ is the probability of the observed sequences
under the assumption that the  underlying process $\{\e_t\}_{t \geq 1}$ passes
through the  point $(n,m)$.  But the length  of the hidden  trajectory remains
unknown when computing $\Q$.  This gives the formula:
\begin{equation}  \label{Qnm}
\Q(X_{1:n}, Y_{1:m})  = \sum_{e  \in \E_{n,m}} \pr  (\e_{1:|e|} =  e, X_{1:n},
Y_{1:m}) .
\end{equation}
Let  us stress that we have
$$
w_t(\theta) =\log \pr(\exists s\geq  1, Z_s =(N_t,M_t); X_{1:N_t}, Y_{1:M_t}),
\quad t\geq 1,
$$
meaning that  the length of the trajectory  is not necessarily $t$,  but is in
fact unknown.\\

$\Q$  is  the  quantity  that  is  computed  by   forward  algorithm  (see
\cite{Durbin}) and which is used  as likelihood in biological applications. It
is computed via recursive equations similar to those of Lemma~\ref{lem:prec}.
In practice, paths  with highest scores according to  the the Needleman-Wunsch
scoring scheme exactly correspond to  highest probability paths in a pair-HMM,
with a corresponding choice of the parameters (\cite{Durbin}).  Thus, the 
quantity  $\Q$   is  used   for  finding  the   best  alignment   between  two
sequences.  Moreover, as  we explained  it in  the introduction,  the  idea of
maximizing this quantity with respect to the parameter $\theta$ has now widely
spread                           among                           practitioners
(\cite{TKF,TKF2,HeinWiuf,Metzler,KnudsenMiyamoto,Miklos}).  The   goal  is  to
obtain an objective choice of  the parameters appearing in the scoring scheme,
taking evolution  into account. Thus, asymptotic properties  of criterion $\Q$
and consequences on  asymptotic properties of the estimator  derived from $\Q$
are of primarily interest. \\
According to  the relation
\eqref{longueur}, asymptotic results  for $t \to  \infty$ will imply
equivalent ones for $n,m \to \infty$. In other words, consistency results obtained
when $t \to \infty$ can be interpreted as valid for long enough observed sequences,
even if one does not know $t$.

\subsection{Biologically motivated restrictions}
Evolution  models  are  commonly  chosen  time reversible,  in  the
limit  of infinitely long  sequences. The reversibility property
implies  that the joint probability of sequence $X$ and an  ancestor
sequence $U$ is not influenced by the fact  that $X$  is a
descendant of sequence  $U$: this  joint probability would be the
same  if $X$ were an ancestor of $U$  or if both were descendants of
a third sequence.  Note that this assumption does not apply on the level of
alignments.  Indeed,   for  single  alignments,  one  may   have  $\pr(\e  =e,
X,Y) \neq \pr(\e=e', Y,X)$, where $e$ and $e'$ are equal on diagonal steps and
have  switched  insertions  and  deletions (namely,  corresponding  paths  are
symmetric around  the axis $x=y$). In fact,  it is the probability  of a whole
given  set of  evolution  events (namely  mutations,  insertions or  deletions
occurring in the  evolution process), which is a  sum over different alignments
$e$ (all representing this same set of evolution events) of probabilities $\pr(\e=e,
X,Y)$,   which   is   conserved   if   we   interchange   the   two   observed
sequences.  More precisely, we always have $\sum_{e \in \E_1} \pr(\e=e, X, Y) =
\sum_{e\in  \E_2} \pr(\e=e,  Y, X)$  where  $\E_1$ and  $\E_2$ are  alignments
subsets representing the same set of evolution events.
\\
Evolution  models  rely  on  two separate  processes:  the  insertion-deletion
(indel)  and the substitution  process and  both are  supposed to  be time
reversible.   As a  consequence of  time reversibility  of indel  process, the
stationary probability of  appearance of an insertion or of  a deletion is the
same, meaning that $p=q$. We thus introduce the following assumption on the
stationary distribution of the hidden Markov chain:
\begin{assumption}
$p=q$.
\end{assumption}
Time  reversibility assumption  on the  substitution process  implies equality
between the marginals of $h$  and individual distributions of the letters,
namely $h_X=f$ and $h_Y=g$. We thus also introduce the following assumption on the emission
distributions:
\begin{assumption}
$h_X=f$ and $h_Y=g$.
\end{assumption}
This last assumption has an interesting consequence on the distribution of only one
sequence:
\begin{lemma}
\label{lemsym}
Under Assumption  3, for any integers  $n$ and $m$, any  $x_{1:n}$ and
any $y_{1:m}$
$$
\pr\left(Z_{t}=(n,m), X_{1:n}=x_{1:n}\right) =
\pr\left(Z_{t}=(n,m)) f^{\otimes n} (x_{1:n}\right),
$$
$$
\pr\left(Z_{t}=(n,m), Y_{1:m}=y_{1:m}\right) =
\pr\left(Z_{t}=(n,m)) g^{\otimes m} (y_{1:m}\right).
$$
\end{lemma}
Here, $f^{\otimes n} (x_{1:n})\triangleq f(x_1) \ldots f(x_n)$.\\

\begin{proof}\\
One has
\begin{eqnarray*}
\lefteqn{ \pr\left(Z_{t}=(n,m), X_{1:n}=x_{1:n}\right)}&& \\
&=&\sum_{y_{1:m}}\pr\left(Z_{t}=(n,m), X_{1:n}=x_{1:n},Y_{1:m}=y_{1:m}\right)\\
&=&
\sum_{e \in \E_{n,m},|e|=t}\sum_{y_{1:m}}\pr\left(\e_{1:t}=e, X_{1:n}=x_{1:n},Y_{1:m}=y_{1:m}\right)\\
&=&
\sum_{e \in \E_{n,m},|e|=t}\pr\left(\e_{1:t}=e\right)
\sum_{y_{1:m}}\pr\left(X_{1:n}=x_{1:n},Y_{1:m}=y_{1:m}|\e_{1:t}=e\right),
\end{eqnarray*}
so that use of equation (\ref{conditional}) and Assumption 3 gives the first assertion of the Lemma.
Proof of the second assertion is similar.
\end{proof}

\section{Information divergence rates}\label{sec:divergence}
\subsection{Definition of Information divergence rates.}
In this section, we investigate the asymptotic properties of the {\it
log-likelihoods}   $\ell_t  (\theta)$   and  $w_t   (\theta)$   when  properly
normalized. We first prove that limiting functions exist. We shall need the following parameter sets $\Theta_{0}$ and $\Theta_{\delta}$, $\delta>0$:
\begin{eqnarray*}
\Theta_{\delta}& =& \left\{\theta \in \Theta \, | \, \pi(i,j) \geq \delta,\,\,
 f(x)\geq \delta,\,\, g(y)\geq \delta, \,\, h(x, y)\geq \delta ,   \forall
 i,j \in \E , \,\,\forall x,y \in {\cal A}\right\},\\
\Theta_{0}&=&\cap_{\delta >0}\Theta_{\delta}
\\
&=&\left\{\theta \in \Theta \, | \, \pi(i,j) >0,\,\,
 f(x)>0,\,\, g(y)>0, \,\, h(x,y)>0 , \,\, \forall i,j \in \E,\,\,\forall x,y
\in {\cal A}\right\}.
\end{eqnarray*}
We shall always assume that $\theta_{0}\in \Theta_{0}$.
\begin{thm}
\label{thdivergence}
The following holds for any $\theta\in\Theta_0$:
\begin{itemize}
\item[i)]
$ t^{-1} \ell_{t} (\theta)$ converges $\pro$-almost surely and in
$\mathbb{L}_1$, as $t$ tends to infinity to
$$
\ell (\theta ) = \lim_{t\rightarrow \infty}\frac{1}{t}\espo \left(\log
\pr(X_{1:N_{t}},Y_{1:M_{t}}) \right) = \sup_{t}\frac{1}{t}\espo
\left(\log \pr(X_{1:N_{t}},Y_{1:M_{t}}) \right).
$$
\item[ii)]
$ t^{-1} w_{t} (\theta)$ converges $\pro$-almost surely and in
$\mathbb{L}_1$, as $t$ tends to infinity to
$$
w (\theta ) = \lim_{t\rightarrow \infty}\frac{1}{t}\espo \left(\log
\Q(X_{1:N_{t}},Y_{1:M_{t}}) \right) = \sup_{t}\frac{1}{t}\espo
\left(\log \Q(X_{1:N_{t}},Y_{1:M_{t}} ) \right).
$$
\end{itemize}
\end{thm}
We then define Information divergence rates:
\begin{defi} $\forall \theta \in \Theta_0,$
$$D (\theta \vert \theta_0)=w (\theta_0)-w (\theta) \quad \text{and} \quad
D^{*}(\theta \vert \theta_0)=\ell (\theta_0)-\ell (\theta).$$
\end{defi}
Note that $D^{*}$ is what is usually called the Information divergence rate in Information Theory:
it is the limit of the normalized Kullback-Leibler divergence between the distributions of
the observations at the true parameter value and another parameter value. However, we also call
$D$ an  Information divergence rate  since $\Q$ may   be interpreted  as a
likelihood.\\

\begin{proof} {\bf of Theorem \ref{thdivergence}}
\\
This proof follows the lines of Leroux (\citeyear{Leroux}, Theorem 2).
We shall use the following version of the sub-additive ergodic
Theorem due to Kingman \citeyear{Kingman} to prove point {\it i)}.
A similar proof may be written for {\it ii)} and is left to the reader.\\
Let $(W_{s,t})_{0\leq s <t}$
be a sequence of random variables such that
\begin{enumerate}
\item
For all $s<t$, $W_{0,t}\geq W_{0,s}+ W_{s,t}$,
\item
For all $k>0$, the joint distributions of $(W_{s+k,t+k})_{0\leq s
<t}$ are the same as those of $(W_{s,t})_{0\leq s <t}$,
\item
$\espo(W_{0,1}^{-}) > -\infty$.
\end{enumerate}
Then $\lim_{t\rightarrow \infty} t^{-1} W_{0,t}$ exists almost
surely. If moreover the sequences $(W_{s+k,t+k})_{k>0}$ are ergodic,
then the limit is almost surely deterministic and equals
$\sup_{t} t^{-1} \mathbb{E}_0(W_{0,t})$. If moreover
$\espo(W_{0,t})\leq At$, for some constant $A$,
then the convergence holds in $\mathbb{L}_1$.\\
We apply this theorem to the auxiliary process
$$
W_{s,t}=   \max_{e    \in   \E}\log   \pr(X_{N_{s}+1:N_{t}},Y_{M_{s}+1:M_{t}}|
\e_{s+1}=e ) +\log(\delta_\theta), \quad 0\leq s<t,
$$
where $\delta_\theta =\min_{e,e'\in \E} \pi(e,e') >0$. 
We are interested in the behaviour of 
$$
U_{s,t}=\log \pr(X_{N_{s}+1:N_{t}},Y_{M_{s}+1:M_{t}}), \quad 0\leq
s<t.
$$
Since we have $ \exp (U_{s,t}) =\sum_{e \in \E} \pr(\e_{s+1} 
=e)\pr(X_{N_s+1:N_t},  Y_{M_s+1:M_t} |\e_{s+1}=e)$  leading to  $  \exp (W_{s,t}
-\log \delta_\theta) \min_{e \in \E}\pr(\e_1=e)
\leq \exp (U_{s,t}) \leq \exp (W_{s,t} -\log \delta_\theta) $, we can
conclude that the desired results on $ \lim_{t \to \infty} t^{-1} U_{0,t}$ and
$ \lim_{t \to  \infty} t^{-1} \espo (U_{0,t})$ follow  from corresponding ones
on the process $W$.

Note that since  $Z_0=(0,0)$ is deterministic, we have  $W_{0,t} = \max_{e \in
  \E}\log \pr (X_{1:N_t}, Y_{1:M_t}|$ $ \e_1=e ) +\log \delta_\theta$. 
Super-additivity (namely point 1.) follows since for any $0\leq
s<t$,
\begin{multline*}
 \pr(X_{1:N_{t}}, Y_{1:M_t}|\e_1=e_1) = \sum_{\substack{e \in \E_{N_t,M_t}
\\|e| =t}} \pr(\e_{2:t} =e_{2:t}, X_{1:N_t}, Y_{1:M_t}|\e_1=e_1)\\
\shoveleft{\geq  \sum_{\substack{e^1 \in \E_{N_s,M_s}\\|e^1|
=s}} \sum_{\substack{e^2 \in \E_{N_{t}-N_s,M_{t}-M_s}\\|e^2| =t-s}}
\pr(\e_{2:s} =e^1_{2:s}, \e_{s+1:t} =e^2, X_{1:N_t}, Y_{1:M_t}|\e_1=e_1)}\\
\shoveleft{ =\sum_{\substack{e^1 \in \E_{N_s,M_s}\\|e^1|
=s}}  \sum_{\substack{e^2 \in \E_{N_{t}-N_s,M_{t}-M_s} \\|e^2| =t-s}}
 \pr(\e_{s+2:t} =e^2_{2:t-s}},
X_{N_{s}+1:N_t}, Y_{{M_s}+1:M_t} | \e_{s+1}=e^2_{1}  )  \\
\shoveright{   \times   \pi(e_s,e_{s+1})  \pr(\e_{2:s}   =e^1_{2:s},X_{1:N_s},
  Y_{1:M_s}|\e_1=e_1 ) }\\ 
=  \sum_{e_s,e_{s+1}\in  \E}\pr(X_{N_{s}+1:N_t}, Y_{{M_s}+1:M_t}  |
\e_{s+1}=e_{s+1} ) \pi(e_s,e_{s+1}) \pr(\e_s=e_s, X_{1:N_{s}},Y_{1:M_s} |\e_1=e_1)\\
\geq \{ \max_{e'\in \E} \pr(X_{N_s+1:N_{t}}, Y_{M_s+1:M_t}|\e_{s+1}=e') \} 
\{ \min_{e,e'}\pi(e,e') \} \pr(X_{1:N_{s}},Y_{1:M_s} |\e_1=e_1) \;, 
\end{multline*}
so that we get $ W_{0,t} \geq W_{0,s} +W_{s,t}$, for any $0\leq s<t$.\\
To understand the distribution of $(W_{s,t})_{0\leq s <t}$, note
that $W_{s,t}$ only depends on trajectories of the random walk going
from the point $(N_{s},M_{s})$  to the point  $(N_{t},M_{t})$ with
length  $t-s$. Since the process $(\e_t)_{t\in \N}$ is
stationary, one gets that the distribution of $(W_{s,t})$ is the
same as that of
$(W_{s+k,t+k})$ for any $k$, so that point $2.$ holds.\\
Point $3.$ comes from:
\begin{eqnarray*}
\espo(W_{0,1}^{-})-\log \delta_\theta &=& \espo \max \{ \log  f(X_1) ; \log 
g(Y_1) ; \log h(X_1,Y_1) \}  > -\infty ,
\end{eqnarray*}
$\pro$-almost surely, since $\theta \in \Theta_0$. Let us fix $0\leq s <t$. The proof that $W^{s,t}=(W_{s+k,t+k})_{k>0}$
is ergodic is the same as that of Leroux (\citeyear{Leroux}, Lemma 1). Let $T$ be
the shift operator, so that if $u=(u_k)_{k\geq 0}$, the sequence
$Tu$ is defined by $(Tu)_{k}=(u)_{k+1}$ for any $k\geq 0$. Let $B$
be an event which is $T$-invariant. We need to prove that $\pro(W^{s,t}
\in B)$ equals $0$ or $1$. For any integer $n$, there exists a
cylinder set $B_n$, depending only on the coordinates $u_k$ with
$-m_n \leq k \leq m_n$ for some sub-sequence $m_n$, such that
$\pro(W^{s,t} \in B\Delta B_{m_n})\leq 1/2^n$. Here, $\Delta$ denotes
the symmetric difference between sets. Since $W^{s,t}$ is stationary
and $B$ is $T$-invariant:
\begin{eqnarray*}
\pro\left(W^{s,t} \in B\Delta B_{m_n}\right)&=&\pro\left(T^{2m_n}W^{s,t} \in B\Delta B_{m_n}\right)\\
&=&\pro\left(W^{s,t} \in B\Delta T^{-2m_n}B_{m_n}\right).
\end{eqnarray*}
Let $\tilde{B}=\cap_{n\geq 1}\cup_{j\geq n}T^{-2m_j}B_{m_j} $.
Borel-Cantelli's Lemma leads to $\pro(W^{s,t} \in B\Delta
\tilde{B})=0$, so that $\pro(W^{s,t}  \in  B)=\pro(W^{s,t}  \in
\tilde{B})=\pro(W^{s,t}  \in  B  \cap  \tilde{B})$.  Now, conditional
on $(\e_t)_{t\in\N}$, the random variables $(W_{s+k,t+k})_{k>
0}$ are strongly mixing, so that the $0-1$  law implies (see
\cite{such})  that for any
fixed sequence  $e$ with values in $\E_\infty$,
the    probability   $\pro(W^{s,t}  \in \tilde{B}\vert
(\e_t)_t=e)$ equals $0$ or $1$, so that
$$
\pro\left(W^{s,t} \in \tilde{B}\right)=P\left((\e_t)_t \in
C\right)
$$
where $C$ is the set of sequences $e$ such that $P(W^{s,t} \in
\tilde{B}\vert (\e_t)_t=e)=1$. But it is easy to see that $C$
is $T$-invariant. Indeed, if $e\in C$ then, since $W^{s,t}$ is
stationary and $\tilde{B}$ invariant,
$$1=\pro(W^{s,t} \in
\tilde{B}\vert (\e_t)_t=e)=\pro(TW^{s,t} \in \tilde{B}\vert
(\e_t)_t=Te)=\pro(W^{s,t} \in \tilde{B}\vert (\e_t)_t=Te)
$$
so that $Te\in C$. Now, since a stationary irreducible Markov chain
is ergodic, $\pro\left((\e_t)_t \in C\right)$ equals $0$ or
$1$. This concludes the proof of ergodicity of the sequence $W^{s,t}$.\\
To end with, note that for any $t\geq 0$, the random variable
$W_{0,t}$ is non positive.
\end{proof}

\subsection{Divergence properties of Information divergence rates}
Information divergence rates should be non negative: this is proved below. They also
should be  positive for parameters  that are different  than the true  one: we
only prove it for some subsets of the parameter set. We thus define
$\Theta_{exp}$  as the  subset of  $\Theta_0$  such that  the expectations  of
$\e_{1}$ under $\theta$ and under $\theta_0$ are not aligned with $(0,0)$:
$$
\Theta_{exp}=\left\{\theta\in\Theta_0\;:\;\forall \lambda >0, \esp (\e_1) \neq \lambda \espo(\e_1)
\right\}.
$$
$\Theta_{marg}$ is the subset of $\Theta_0$ such that Assumption 3 holds:
$$
\Theta_{marg}=\left\{\theta\in\Theta_0\;:\;h_{X}=f,\;h_{Y}=g
\right\}.
$$
\begin{thm}
\label{contrast} Information divergence rates satisfy:
\begin{itemize}
\item
For all $\theta \in \Theta_0$, $D(\theta \vert \theta_0) \geq 0$ and $D^{*} (\theta \vert \theta_0) \geq 0$.
\item
For any $\theta  \in \Theta_{exp}$, $\theta \neq \theta_0$, we have $D(\theta
\vert \theta_0) > 0$ and $D^{*}(\theta \vert \theta_0) > 0$.
\item
If $\theta_0$ and $\theta$ are in $\Theta_{marg}$,
$D(\theta
\vert \theta_0) > 0$ and $D^{*}(\theta \vert \theta_0) > 0$
as soon as $f\neq f_0$ or $g\neq g_0$.
\end{itemize}
\end{thm}
Notice that in case Assumption 2 holds, the expectations of $\e_{1}$ under
$\theta$ and under $\theta_0$ are aligned with $(0,0)$. In this case, we were not able
to prove that $h\neq h_0$ implies positivity of information divergence rates.\\
$\\$
\begin{proof}\\
Since for all $t$,
$$
\espo \left(\log
\pro(X_{1:N_{t}},Y_{1:M_{t}}) \right)-
\espo \left(\log
\pr(X_{1:N_{t}},Y_{1:M_{t}}) \right)
$$
is a Kullback-Leibler divergence, it is non negative, and the limit
$D^{*}(\theta \vert \theta_0)$ is also non negative.\\
Let us prove that $D(\theta \vert \theta_0)$ is also  non negative.
To compute  the value  of the expectation  $\espo [w_{t}  (\theta )]$,
introduce the set $A_{t}$ of all possible values of $Z_t$:
\begin{equation}
\label{Zset}
A_{t}=\left\{(n,m)\in \N ^{2}\;:\;n\vee m \leq t \leq n+m
\right\}.
\end{equation}
Then,
\begin{equation}
\label{ltheta}
\espo [w_{t} (\theta )]=\sum_{(n,m)\in A_t} \sum_{x_{1:n},y_{1:m}}
\pro(Z_{t} =(n,m), X_{1:n}=x_{1:n},Y_{1:m}=y_{1:m}  )
\log \Q(x_{1:n},y_{1:m}).
\end{equation}
Now, by definition,
$$
D\left(\theta \vert \theta_0 \right )
= \lim_{t\rightarrow +\infty} \frac{1}{t} \espo
\left(\log \frac{\Qo(X_{1:N_{t}},Y_{1:M_{t}} )}{\Q(X_{1:N_{t}},Y_{1:M_{t}} )}\right).
$$
By using Jensen's inequality,
\begin{equation*}
\espo
\left(\log \frac{\Q(X_{1:N_{t}},Y_{1:M_{t}})}{\Qo(X_{1:N_{t}},Y_{1:M_{t}})}\right)
\leq
\log \espo \left(\frac{\Q(X_{1:N_{t}},Y_{1:M_{t}})}{\Qo(X_{1:N_{t}},Y_{1:M_{t}})}
\right).
\end{equation*}
But
\begin{eqnarray*}
&& \espo\left(\frac{\Q(X_{1:N_{t}},Y_{1:M_{t}})}{\Qo(X_{1:N_{t}},Y_{1:M_{t}})}
\right)\\
&=& \sum_{(n,m)\in A_t} \sum_{x_{1:n},y_{1:m}}
\pro(Z_{t} =(n,m), X_{1:n}=x_{1:n},Y_{1:m}=y_{1:m}  )
\times \frac{\Q(x_{1:n},y_{1:m})}{\Qo(x_{1:n},y_{1:m})}\\
&\stackrel{(a)}{\leq} & \sum_{(n,m)\in A_t} \sum_{x_{1:n},y_{1:m}} \pr(\exists
s \geq 1, Z_s =(n,m), X_{1:n}=x_{1:n}, Y_{1:m}=y_{1:m})\\
&\leq &\sum_{(n,m)\in A_t} \pr(\exists s\geq 1, Z_{s} =(n,m)),
\end{eqnarray*}
where $(a)$ comes from expression (\ref{Qnm}). Finally,
\begin{equation}
\label{majobrute}
\lim_{t\rightarrow +\infty}\frac{1}{t}\left(w_{t}(\theta)-w_{t}(\theta_0)
\right)\leq \liminf_{t\rightarrow +\infty}\frac{1}{t}
\log\left[\sum_{(n,m)\in A_t} \pr \left(\exists s\geq 1, Z_{s} =(n,m) \right)
\right] .
\end{equation}
But the cardinality of $A_t$ is at most $t^2$, so that
$$
\lim_{t\rightarrow +\infty}\frac{1}{t}\left(w_{t}(\theta)-w_{t}(\theta_0)
\right)\leq \liminf_{t\rightarrow +\infty}\frac{1}{t}
\log t^2 =0,
$$
and
$$
\forall \theta \in \Theta_0, \;D(\theta \vert \theta_0) \geq 0.
$$

Since $\theta\in \Theta_0$, there exists $\delta_{\theta}$ such that
$\theta\in \Theta_{\delta_{\theta}}$.
By using (\ref{Qnm}), one gets the lower bound
$$
\Q(x_{1:n},y_{1:m} ) \geq
\delta_{\theta}^{n + m}\inf_{e\in \E_{n,m}} \left[ \pr\left(\e_{1:|e|} =e
  \right)\right].
$$
Since trajectories $e$ in $\E_{n,m}$ have length at most $n+m$,
$$
\inf_{e\in \E_{n,m}} \left [ \pr\left(\e_{1:|e|} =e \right)\right]
\geq \delta_{\theta}^{n + m}.
$$
Note also  that if $(n,m)$  belongs to  $A_t$ then we  have $n+m \leq  2t$ and
$n\vee m \geq t/2$.
Thus, uniformly with  respect  to $(n,m)\in  A_t$  and  to $x_{1:n}$  and
$y_{1:m}$,
\begin{equation}
\label{Qcoince}
4t \log \delta_{\theta} \leq \log \Q(x_{1:n},y_{1:m} ) \leq 0 .
\end{equation}
Moreover, with
$$
\rho_{\theta} = \|f\|_{\infty} \vee \|g\|_{\infty} \vee \|h\|_{\infty}
\leq 1- \delta_{\theta}<1
$$
one has for any integers $n$, $m$, any $x_{1:n}$ and $y_{1:m}$
$$
\Q(x_{1:n},y_{1:m} ) \leq \rho_{\theta}^{n\vee m} .
$$
In this case,  for all $t$,
and uniformly with respect to $(n,m)\in A_t$ and to $x_{1:n}$ and $y_{1:m}$,
\begin{equation}
\label{Qcoince2}
\log \Q(x_{1:n},y_{1:m} ) \leq  \frac{t}{2}\log (1- \delta_{\theta}).
\end{equation}
Inequalities \eqref{Qcoince} and \eqref{Qcoince2} allow to conclude that
$$
-C_{\theta_0} \leq w (\theta_0 ) \leq -c_{\theta_0} \quad\text{and} \quad -C_{\theta} \leq w (\theta ) \leq -c_{\theta} .
$$
Then, as soon as $B_t$ is a
set such that
\begin{equation}
\label{bt}
\lim_{t\rightarrow +\infty} \pro\left(Z_{t} \notin B_{t}  \right) =0,
\end{equation}
we have
$$
D\left(\theta \vert \theta_0 \right )
= \lim_{t\rightarrow +\infty} \frac{1}{t} \espo \left[
\left(\log \frac{\Qo(X_{1:N_{t}},Y_{1:M_{t}} )}{\Q(X_{1:N_{t}},Y_{1:M_{t}} )}\right)
\1{Z_{t} \in B_{t}}\right].
$$
Now, using Jensen's inequality,
\begin{equation*}
\espo \left[
\left(\log \frac{\Q(X_{1:N_{t}},Y_{1:M_{t}})}{ \Qo(X_{1:N_{t}},Y_{1:M_{t}})} \right)
\1{Z_{t} \in B_{t}}\right]
\leq \pro \left(Z_{t} \in B_{t} \right)
\log         \espo        \left(\frac{\Q(X_{1:N_{t}},Y_{1:M_{t}})}{        \Qo
    (X_{1:N_{t}},Y_{1:M_{t}})}
\Big\vert Z_{t}\in B_{t} \right).
\end{equation*}
But as previously seen,
\begin{eqnarray*}
&& \espo\left(\frac{\Q(X_{1:N_{t}},Y_{1:M_{t}})}{\Qo(X_{1:N_{t}},Y_{1:M_{t}})}
\vert Z_{t}\in B_{t}\right)\\
&=& \sum_{(n,m)\in B_t} \sum_{x_{1:n},y_{1:m}}
\frac{\pro(Z_{t} =(n,m), X_{1:n}=x_{1:n},Y_{1:m}=y_{1:m}  )}{\pro(Z_{t} \in
  B_{t} )}
\times \frac{\Q(x_{1:n},y_{1:m})}{\Qo(x_{1:n},y_{1:m})}\\
&\leq &\sum_{(n,m)\in B_t} \frac{\pr(\exists s\geq 1, Z_{s} =(n,m))}
{\pro(Z_{t} \in B_{t})} .
\end{eqnarray*}
Finally,
\begin{equation}\label{majo2}
 - D\left(\theta \vert \theta_0 \right )
\leq \lim_{t\rightarrow +\infty} \frac{1}{t} \log \pr(\exists s\geq 1, Z_s \in B_t).
\end{equation}

Let us now consider the case where the expectations of $\e_1$ under
parameters $\theta$ and $\theta_0$ are not aligned with $(0,0)$, that is
$\theta \in \Theta_{exp}$.
We have
$$
\eta = \inf_{\lambda \in \R} \left\|\esp (\e_1)-\lambda
\espo (\e_1)\right\| >0,
$$
where $\|\cdot \|$ denotes the euclidean norm. Define
$$
B_{t}=\left\{(n,m) \in A_t \;:\;\left\|\frac{(n,m)}{t}-
\espo (\e_1)\right\| \leq \frac{\eta}{4}
\right\}.
$$
Then, (\ref{bt}) holds.  Any trajectory $e$ ending at point $(n,m)$
has length at least $n\vee m$ which is at least $t/2$ when $(n,m)\in B_{t}$.
Thus for such $(n,m)$:
\begin{eqnarray*}
\pr \left(\exists s\geq 1, Z_{s} =(n,m) \right)& \leq &
\pr \left(\exists s \geq \frac{t}{2},
\inf_{\lambda \in \R} \left\|\frac{Z_{s}}{s}-\lambda
\espo   (\e_1)\right\|    \leq   \frac{t}{s}   \left\|\frac{Z_{s}}{t}-   \espo
(\e_1)\right\| \right)\\
& \leq &
\pr \left(\exists s \geq \frac{t}{2},
\inf_{\lambda \in \R} \left\|\frac{Z_{s}}{s}-\lambda
\espo (\e_1)\right\| \leq \frac{\eta}{2} \right)\\
&\leq&
\pr\left(\exists s\geq \frac{t}{2},
\left\|\frac{Z_{s}}{s}-
\esp (\e_1)\right\| \geq \frac{\eta}{2} \right).
\end{eqnarray*}
Now, using easy Cramer-Chernoff bounds, since $\pi$ is irreducible,
one has that there exists a positive $c(\eta)$ and some $s_0>0$ such that
as soon as $s\geq s_0$,
$$
\pr \left(
\left\|\frac{Z_{s}}{s}-
\esp (\e_1)\right\| \geq \frac{\eta}{2} \right)
\leq \exp \left(-s c(\eta)\right),
$$
and by summing over $s$, there also  exists a positive $C$ such that for large
enough $t$,
$$
\pr \left( \exists s\geq \frac{t}{2} \;:\;
\left\|\frac{Z_{s}}{s}-
\esp (\e_1)\right\| \geq \frac{\eta}{2} \right)
\leq C\exp \left(-t c(\eta)/2\right).
$$
Thus, using (\ref{majo2}), one obtains that for $\theta \in \Theta_{exp}$:
$$
D(\theta \vert \theta_0) \geq \frac{c(\eta)}{2}>0.
$$

Let us now consider the case where  $\theta_0$ and $\theta$ are in $\Theta_{marg}$.
Then, using Jensen's Inequality and definition \eqref{Qnm},
\begin{eqnarray*}
&&   \espo  \left(  \log
  \frac{\Q(X_{1:N_t} ,Y_{1:M_t}) } {\Qo(X_{1:N_t}, Y_{1:M_t})}\right)  \\
&=&  \sum_{(n,m)  \in A_t}  \sum_{x_{1:n}}
\sum_{y_{1:m}} \pro(Z_t =(n,m), X_{1:n}=x_{1:n}, Y_{1:m}=y_{1:m}) \log
\frac{\Q(x_{1:n} ,y_{1:m}) } {\Qo(x_{1:n} ,y_{1:m})} \\
&\leq & \sum_{(n,m) \in A_t} \sum_{x_{1:n}}
 \pro(Z_t =(n,m), X_{1:n}=x_{1:n}) \\
&& \log \left(\sum_{y_{1:m}}
\frac{\pro(Z_t    =(n,m),    X_{1:n}=x_{1:n},   Y_{1:m}=y_{1:m})\Q(x_{1:n} ,y_{1:m}
) }  {\pro(Z_t =(n,m), X_{1:n}=x_{1:n})  \Qo(x_{1:n} ,y_{1:m}) }
\right) \\
&\leq & \sum_{(n,m) \in A_t} \sum_{x_{1:n}}
 \pro(Z_t =(n,m)) f_0^{\otimes n} (x_{1:n})
 \log \left(\frac{ \pr(\exists s\geq 1, Z_s=(n,m)) f^{\otimes n} (x_{1:n}) } {\pro(Z_t
=(n,m)) f_0^{\otimes n} (x_{1:n}) } \right),
\end{eqnarray*}
where the last inequality comes from Lemma~\ref{lemsym} and the fact that $\pro(Z_t =(n,m),
X_{1:n}=x_{1:n}, Y_{1:m}=y_{1:m}) \leq \Qo(x_{1:n}, y_{1:m})$.\\
Thus, since $t^{-1}N_t$ tends to $(1-p)$, $\pro$-a.s. as $t$ tends to infinity, and $(1-p)>0$ since
$\theta\in\Theta_0$, we have
\begin{multline*}
  -D(\theta|\theta_0)   \leq       \limsup_{t   \to   +\infty}   \frac{1}{t}
  \sum_{(n,m)\in A_t,n\geq \frac{(1-p)}{2}t} \pro(Z_t=(n,m)) \Big\{\log \frac {\pr(\exists s\geq 1,
      Z_s=(n,m))} {\pro(Z_t =(n,m)) } \\
+ \frac{(1-p)}{2}t \sum_x f_0(x)\log \frac{f(x)}{f_0(x)}
\Big\}
\leq \frac{(1-p)}{2}\sum_x f_0(x)\log \frac{f(x)}{f_0(x)} <0,
\end{multline*}
as soon as $f\neq f_0$.
A similar proof applies if $g\neq g_0$.\\
Proofs of divergence properties for $D^{*}$ follow the same lines.
\end{proof}

\subsection{Continuity properties}

On $\Theta_{\delta}$, the log-likelihoods are uniformly equicontinuous, with a modulus of continuity
that does not depend on trajectories, as appears in the proof of the following Lemma.
\begin{lemma}
\label{lemequi}
The families of functions $\{ t^{-1} w_{t} (\theta )\}_{t \geq 1}$ and
$\{ t^{-1} \ell_{t} (\theta )\}_{t \geq 1}$ are uniformly
equicontinuous on $\Theta_{\delta}$.
\end{lemma}
A consequence of this Lemma and the compactness of $\Theta_{\delta}$ is:
\begin{cor} \label{corconti} The following holds:
\begin{itemize}
\item [i)] $\{t ^{-1} w_t(\theta)\}_t$ (resp. $\{ t^{-1} \ell_t(\theta)\}_t$) converges $\pro$-almost surely to
$w(\theta)$ (resp. to $\ell (\theta)$) uniformly on $\Theta_{\delta}$;
\item [ii)] $\ell (\theta)$ and $w(\theta)$ are uniformly continuous on $\Theta_{\delta}$.
\end{itemize}
\end{cor}

\begin{proof} {\bf of Lemma \ref{lemequi}}\\
Let $ \alpha>0$,
and $\theta_1, \theta_2 \in \Theta_{\delta}$ such that  $\|\theta_1
-\theta_2\|_{\infty}\leq \alpha$.\\
Let us denote $\mu_{\theta_i}$, $\pi_{\theta_i}$, $f_{\theta_i}$,
$g_{\theta_i}$ and
$h_{\theta_i}$  the parameters of the hidden Markov chain and of
the emission distributions under $\theta_i$, $i=1,2$. \\
For any $e \in \E_{N_t,M_t}$ :
\begin{multline*}
\frac{1}{t} \left|\log \pru (\e_{1:|e|}=e,X_{1:N_t},Y_{1:M_t})-\log
\prd(\e_{1:|e|}=e,X_{1:N_t},Y_{1:M_t}) \right|\\
\shoveleft{ \leq \frac{1}{t} \left|\log
\mu_{\theta_1}(e_1)\!-\!\log \mu_{\theta_2}(e_1)\right| \!+
\frac{1}{t}\sum_{k,l \in \E}
\!\Big( \sum_{i=2}^{|e|}\!\1{e_{i-1}\!\!=\!k,e_i\!\!=\!l}\!\Big)
\left|\log \pi_{\theta_1}(k,l)\!-\!\log \pi_{\theta_2}(k,l)\right| } \\
\shoveleft{ +\frac{1}{t}\sum_{a\in{\cal A}} \left\{ \Big( \sum_{i=1}^{|e|}
\1{e_i=(1,0),X_{N_i}=a}\Big) \,|\log
f_{\theta_1}(a)-\log f_{\theta_2}(a)| \right.}\\
\shoveright{\left.+ \Big( \sum_{i=1}^{|e|}
\1{e_i=(0,1),Y_{M_i}=a}\Big) \,|\log
\,g_{\theta_1}(a)-\log \,g_{\theta_2}(a)|\right\}}  \\
+\frac{1}{t}\sum_{a,a'\in{\cal A}} \Big( \sum_{i=1}^{|e|}
\1{e_i=(1,1),X_{N_i}=a,Y_{M_i}=a'}\Big) \,\,|\log
h_{\theta_1}(a,a')-\log h_{\theta_2}(a,a')|.
\end{multline*}
In this sum, at most $2|e|$ terms are non null.
Since all the components of $\theta_i$, $i=1,2$ are
bounded below by $\delta$ and $\|\theta_1-\theta_2\|_\infty\leq \alpha$, we have :
$$
\frac{1}{t}|\log \pru(\e_{1:|e|}=e,X_{1:N_t},Y_{1:M_t})-\log
\prd(\e_{1:|e|}=e,X_{1:N_t},Y_{1:M_t})| \leq \frac{2|e|}{t}\frac{\alpha}{\delta}.
$$
But for any $e \in \E_{N_t,M_t}$, we have $|e|\leq 2t$, so that
$$\frac{1}{t}|\log \pru(\e_{1:|e|}=e,X_{1:N_t},Y_{1:M_t})-\log
\prd(\e_{1:|e|}=e,X_{1:N_t},Y_{1:M_t})| \leq \frac{4\alpha}{\delta}, $$
as soon as $\|\theta_1
-\theta_2\|_{\infty}\leq \alpha $.\\
Now we get
\begin{eqnarray*}
\Qu(X_{1:N_t},Y_{1:M_t}) & = & \sum_{e\in \E_{N_t,M_t}}
\pru(\e_{1:|e|}=e,X_{1:N_t},Y_{1:M_t}) \\
& & \leq \exp\left\{\frac{4\alpha}{\delta}t \right\} \sum_{e\in \E_{N_t,M_t}}
\prd(\e_{1:|e|}=e,X_{1:N_t},Y_{1:M_t}) \\
& & \leq \exp \left\{\frac{4\alpha}{\delta}t \right \} \,\,\Qd(X_{1:N_t},Y_{1:M_t}),
\end{eqnarray*}
and $t^{-1} \log \Qu(X_{1:N_t},Y_{1:M_t}) \leq
4\alpha /\delta + t^{-1} \log \Qd(X_{1:N_t},Y_{1:M_t})$. Since
this is symmetric in $\theta_1$ and $\theta_2$, one obtains that for any
$\theta_1, \theta_2 \in \Theta_{\delta}$ such that
$\|\theta_1
-\theta_2\|_{\infty}\leq \alpha $,
$$
|\frac{1}{t}w_{t}(\theta_{1}) - \frac{1}{t}w_{t}(\theta_{2})|\leq \frac{4\alpha}{\delta}.
$$
The same proof applies to $t ^{-1} \ell_{t}$. 
\end{proof}

\section{Statistical properties of estimators}
We now want to focus on a  particular form of the pair-HMM, relying on a
re-parametrization of the  model. Indeed, the pair-HMM  has been introduced
to  take  into account  evolutionary  events.  The corresponding  evolutionary
parameters are  the ones of interest  and practitioners aim  at estimating those
parameters rather than the full pair-HMM.
Examples   of  such   re-parametrization   may  be   found   for  instance   in
\cite{TKF,TKF2} (see also Section 5 of this paper).
Let $\beta \mapsto \theta (\beta)$ be a continuous parametrization from
some set $B$ to $\Theta$. For any $\delta >0$, let $B_{\delta}=\theta^{-1}(\Theta_{\delta})$.
We assume that $\beta_0 =\theta^{-1}(\theta_0)$ in $B_{\delta}$ for some $\delta >0$.
Use of pair-HMM algorithms  to estimate evolutionary parameters corresponds to
the estimator
\begin{defi}
$$
\widehat{\beta}_{t}=\mathop{\textrm{ Argmax }}_{\beta \in B_{\delta}} w_{t}(\theta (\beta)).
$$
\end{defi}
Then,
\begin{thm}
\label{thcons}
If the set of maximizers of $w(\theta (\beta))$ over $B_{\delta}$ reduces to $\{\beta_0\}$,
$\widehat{\beta}_{t}$
converges $\pro$-almost surely to $\beta_0$.
\end{thm}
The proof of this theorem follows from Corollary \ref{corconti} and usual arguments for M-estimators.\\
The  condition  that  the  set  of  maximizers  of  $w(\theta  (\beta))$  over
$B_{\delta}$  reduces  to $\{\beta_0\}$  corresponds  to some  identifiability
condition and thus may not be avoided. \\

Another interesting approach to sequence alignment by pair-HMMs is to
consider a non-informative prior distribution on the parameters to
produce, via a MCMC procedure, the posterior distribution of the alignments
and parameters given the observed sequences. \\
Using $\Q$
as the likelihood of the observed sequences produces a posterior distribution as follows.
Let $\nu$ be a prior probability measure on $B_{\delta}$
and $\bar{\beta}$ a random vector distributed according to
$\nu$. MCMC algorithms approximate the random distribution $\nu_{| X_{1:N_t},Y_{1:M_t}}$
interpreted as the
posterior measure given observations $X_{1:N_t}$ and $Y_{1:M_t}$:
\begin{equation}
\frac{\Qb(X_{1:N_t},Y_{1:M_t})\nu(d\beta)}{\int_{B_{\delta}}\Qb(X_{1:N_t},Y_{1:M_t})\nu(d\beta)}.
\end{equation}
This  leads to  Bayesian consistent  estimation of  $\beta_0$ as  in classical
statistical models (see \cite{IH} for instance). 
Notice that since $w_t$ is
not  the logarithm  of a  probability distribution  on the  observation space,
these results are not direct consequences of classical ones. Though, the proof
follows classical ideas of Bayesian theory.
\begin{thm}
\label{thbayes}
If the set of maximizers of $w(\theta (\beta))$ over $B_{\delta}$ reduces to $\{\beta_0\}$,
and if $\nu$ weights $\beta_0$, then
the sequence of posterior measures $\nu_{| X_{1:N_t},Y_{1:M_t}}$
converges in distribution
$\pro$-almost surely to the Dirac mass at $\beta_0$.
\end{thm}
\begin{proof}
Let $m: B_{\delta} \rightarrow \R$ be any continuous, bounded function.
For any $\epsilon >0$, let $\alpha$ such that $|m(\beta )-m(\beta' )|\leq \epsilon$
as soon as $\|\beta -\beta'\| \leq \alpha$.
We
have
\begin{multline*}
\left|\int_{B_{\delta}}m(\beta)
\frac{\Qb(X_{1:N_t},Y_{1:M_t})}
{\int_{B_{\delta}}\Qb(X_{1:N_t},Y_{1:M_t})\nu(d\beta)}
\nu(d\beta)- m(\beta_0)\right| \\
\leq\frac{\int_{B_{\delta}}
|m(\beta)-m(\beta_0)|
\Qb(X_{1:N_t},Y_{1:M_t})\nu(d\beta)}{\int_{B_{\delta}}
\Qb(X_{1:N_t},Y_{1:M_t})\nu(d\beta)}  .
\end{multline*}
But
$$
\frac{\int_{\|\beta -\beta_0\|\leq \alpha }
|m(\beta)-m(\beta_0)|
\Qb(X_{1:N_t},Y_{1:M_t})\nu(d\beta)}{\int_{B_{\delta}}\Qb(X_{1:N_t},Y_{1:M_t})
\nu(d\beta)}
\leq \epsilon
$$
so that
\begin{eqnarray*}
\lefteqn{\left|\int_{B_{\delta}}m(\beta)
\frac{\Qb(X_{1:N_t},Y_{1:M_t})}
{\int_{B_{\delta}}\Qb(X_{1:N_t},Y_{1:M_t})\nu(d\beta)}
\nu(d\beta)-m(\beta_0)\right|} && \\
&\leq&
\epsilon +
2 \|m\|_{\infty}\frac{\int_{\|\beta -\beta_0 \| > \alpha}
\Qb(X_{1:N_t},Y_{1:M_t})\nu(d\beta)}{\int_{B_{\delta}}\Qb(X_{1:N_t},Y_{1:M_t})
\nu(d\beta)}
\\
&=&\epsilon +2
\|m\|_{\infty}\frac{\int_{\|\beta-\beta_0 \| > \alpha } \exp \left\{t \left(
\frac{1}{t}w_t(\theta (\beta))\right)\right\}
\nu(d\beta)}{\int_{B_{\delta}} \exp\left\{t
\left(\frac{1}{t}w_t(\theta (\beta))\right)\right\}\nu(d\beta)}.
\end{eqnarray*}
Use of Corollary \ref{corconti} and the fact that the set of maximizers of $w(\theta (\beta))$
over $B_{\delta}$ reduces to $\{\beta_0\}$ gives $\eta >0$ and $T$ such that for $t>T$ and
$\|\beta-\beta_0 \| > \alpha$, $t ^{-1} w_t(\theta (\beta))- t^{-1} w_t(\theta (\beta_0))\leq -\eta$,
and then there exists $\gamma >0$ such that for
$t>T$ and
$\|\beta-\beta_0 \| \leq \gamma$, $t ^{-1} w_t(\theta (\beta))- t^{-1} w_t(\theta (\beta_0))\geq -\frac{\eta}{2}$. 
Then
\begin{eqnarray*}
\lefteqn{
\frac{\int_{\|\beta-\beta_0 \| > \alpha } \exp \left\{t \left(
\frac{1}{t}w_t(\theta (\beta))\right)\right\}
\nu(d\beta)}{\int_{B_{\delta}}\exp\left\{t
\left(\frac{1}{t}w_t(\theta(\beta))\right)\right\}\nu(d\beta)}} &&\\
&\leq&
\frac{\int_{\|\beta-\beta_0 \| > \alpha } \exp \left\{t \left(
\frac{1}{t}w_t(\theta(\beta))-\frac{1}{t}w_t(\theta(\beta_0))\right)\right\}
\nu(d\beta)}{\int_{\|\beta-\beta_0 \| \leq \gamma}
\exp\left\{t
\left(\frac{1}{t}w_t(\theta (\beta))-\frac{1}{t}w_t(\theta(\beta_0))\right)\right\}\nu(d\beta)}\\
&\leq& \left(\exp \left\{-t\frac{\eta}{2}\right\}\right)
\frac{\int_{\|\beta-\beta_0 \| > \alpha } \nu(d\beta)}
{\int_{\|\beta-\beta_0 \| \leq \gamma } \nu(d\beta)} .
\end{eqnarray*}
Using that $\nu$ weights $\beta_0$
we finally obtain
\begin{equation}\label{limite}
\lim_{t\rightarrow \infty} \left|\int_{B_{\delta}}m(\beta)
\frac{\Qb(X_{1:N_t},Y_{1:M_t})}{\int_{B_{\delta}}\Qb(X_{1:N_t},Y_{1:M_t})\nu(d\beta)}
\nu(d\beta) - m(\beta_0)\right|=0 \quad \pro-a.s.
\end{equation}
But
it exists a countable collection  of continuous and bounded functions that
are determining for convergence in distribution and the union of the
corresponding null sets in which \eqref{limite} does not hold is
still a null set. Then
\begin{equation}\label{Rconv}
\nu_{ | X_{1:N_t},Y_{1:M_t}} \leadsto \delta_{\beta_0}
\quad \pro-a.s.
\end{equation}
\end{proof}

\section{Simulations}
\subsection{A simple model}
For the whole simulation procedure we consider the following
substitution model:
\begin{equation}
h(x,y) =\left\{
  \begin{array}{ll}
f(x) (1-e^{-\alpha})f(y) & \text{ if } x\neq y\\
f(x)\{ (1-e^{-\alpha})f(x) +e^{-\alpha} \} & \text{ otherwise, }
  \end{array}
\right .
\end{equation}
where $\alpha  >0$ is called the  substitution rate and for  every
letter $x$, $f(x)$ equals the equilibrium probability of $x$. This
equilibrium probability distribution  is assumed to be known  and
will not be  part of the parameter.  Here, the  emission
distribution $g$ equals $f$,  and Assumption 3 holds. The unknown
parameter is thus $\beta=(\pi,\alpha)$. This is a classical
substitution model (used for instance in \cite{TKF}) where the
substitution rate is independent of the type of nucleotide being
replaced and
$1-e^{-\alpha}$ represents the probability that a substitution
occurs. We shall consider hidden Markov chains that satisfy
Assumption 2, and will present:
\begin{itemize}
\item
Simulations with i.i.d. $(\e_s)_s$ where probabilities of horizontal
or vertical moves equal $p_0$ and probability of diagonal moves
equals $r_0=1-2p_0$. Here, the parameter reduces to
$\beta=(p,\alpha)$.
\item
Simulations with stationary Markov chains  such that $p_0 =q_0$. The parameter
dimension then reduces to $6$ (including $\alpha$).
\end{itemize}
Notice that none of these situations is covered by Theorem \ref{contrast}: we do not
know in those cases whether the information divergence rates are positive at a
parameter value different from the true one.\\
In both cases, we get estimations of the parameters via MLE (taking
$\Q$ as the likelihood as it is done in practice), and in the
i.i.d. case we compute and
compare the functions $w$ and $\ell$. 

\subsection{Simulations with i.i.d. $(\e_s)_s$}
\begin{figure}[!h]
 \centering
 \includegraphics[width=8.2cm]{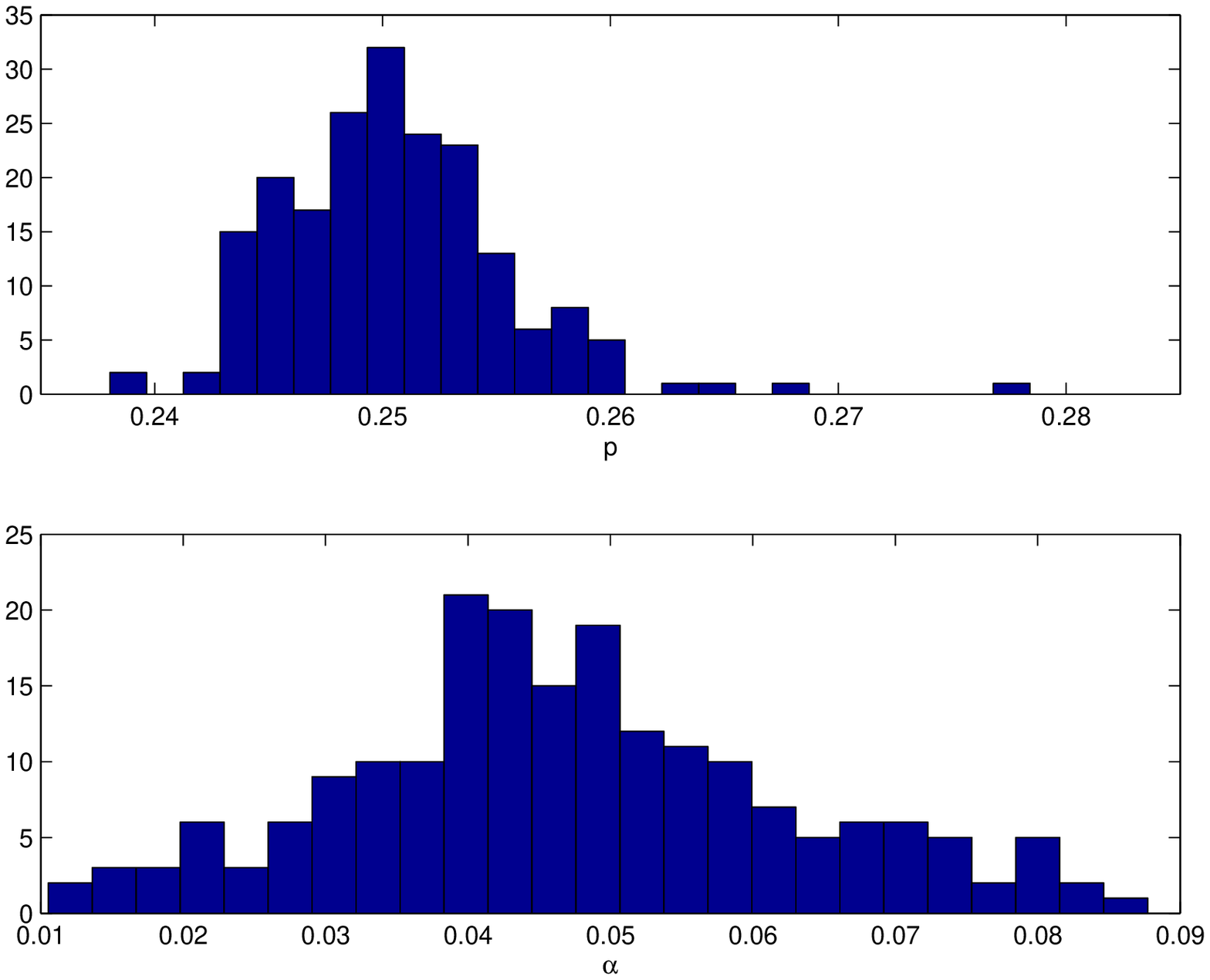}
 \includegraphics[width=6.8cm]{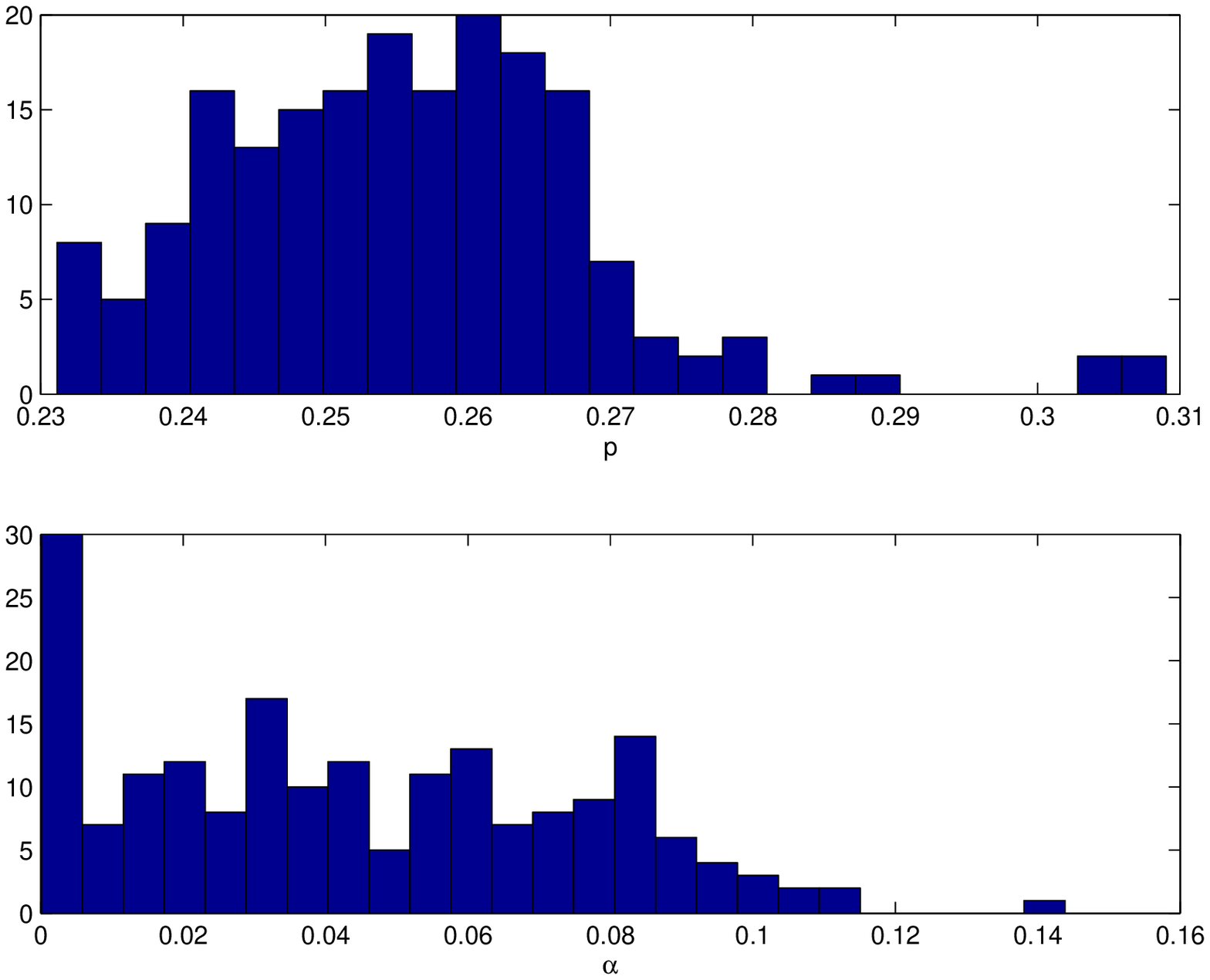}
 \caption{Histograms of maximum likelihood estimations of parameters obtained with  200 simulations
  from the i.i.d. model.
  On the left: estimation of $p$
  given $\alpha=\alpha_0$ and estimation of $\alpha$ given $p=p_0$. On the right: joint estimation of $p$ and $\alpha$.}
\label{iid}
\end{figure}

We have simulated 200 alignments of length 15000 with substitution
rate $\alpha_0=0.05$ and $p_0=q_0=0.25$. We have set the equilibrium
probability of every nucleotide to $0.25$. We show in Figure \ref{iid}
histograms for the maximum likelihood estimations of both
parameters. In a first part we keep $\alpha$  fixed at $\alpha_0$ and
estimate $p$ and then we keep $p$  fixed at $p_0$ and estimate
$\alpha$. That produces good estimations of the parameters even if
$\alpha$ is a bit underestimated. However when estimating $p$ and
$\alpha$ simultaneously (second part) we obtain no satisfying results especially
on $\alpha$ (see Figure \ref{iid}). 

\begin{figure}[!h]
 \centering
 \includegraphics[width=7.6cm]{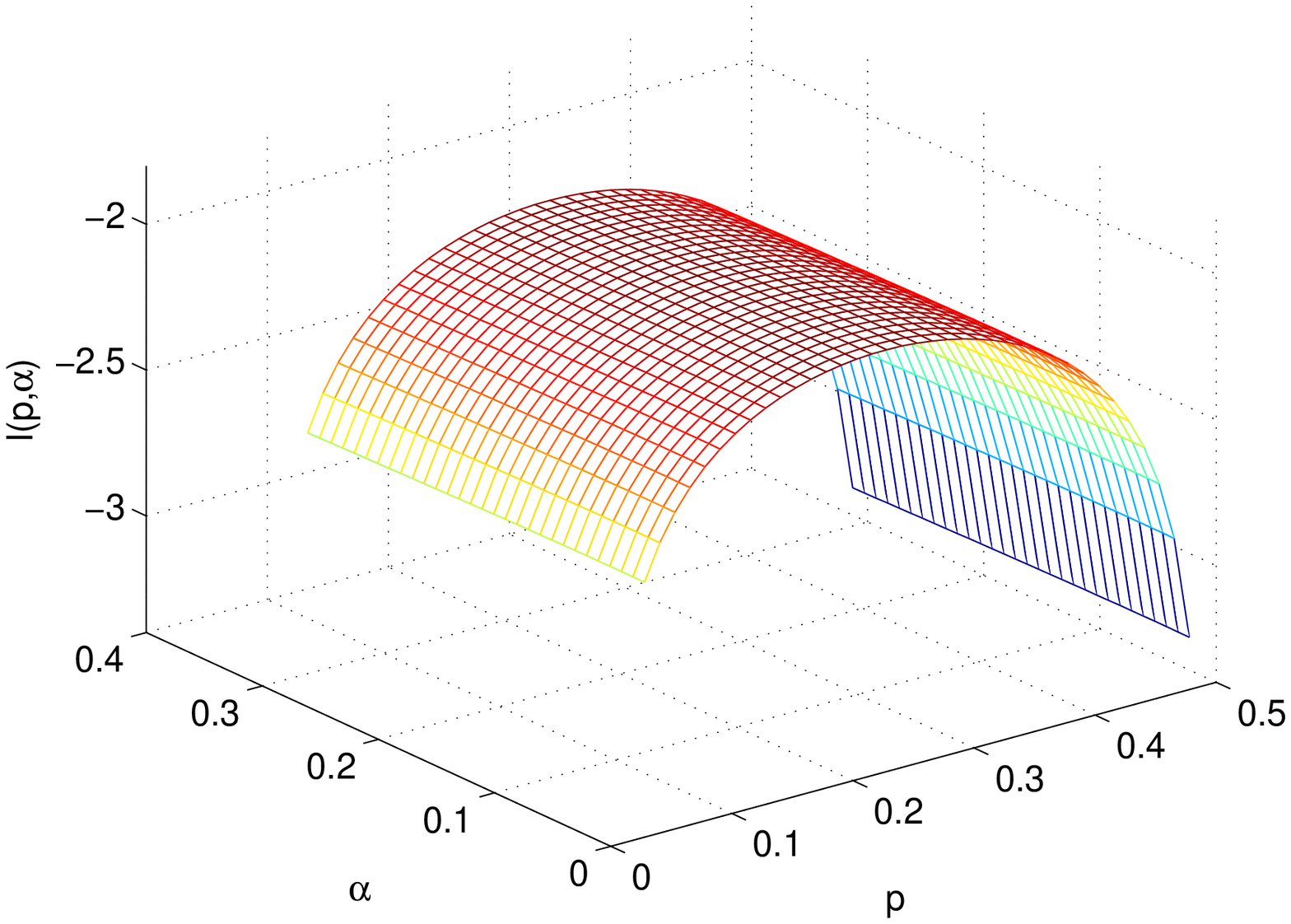}\includegraphics[width=7.6cm]{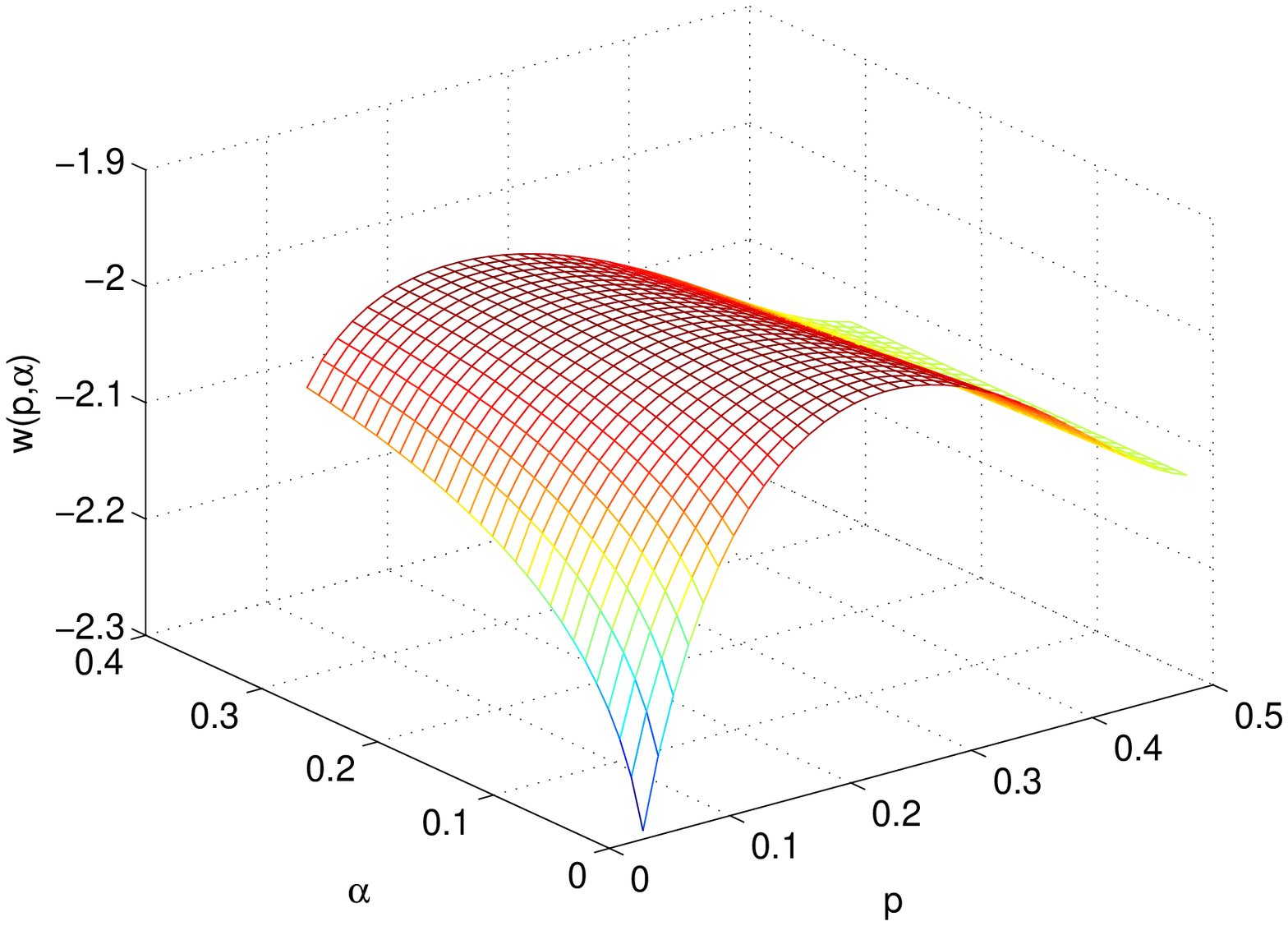}\\
 \includegraphics[width=15cm]{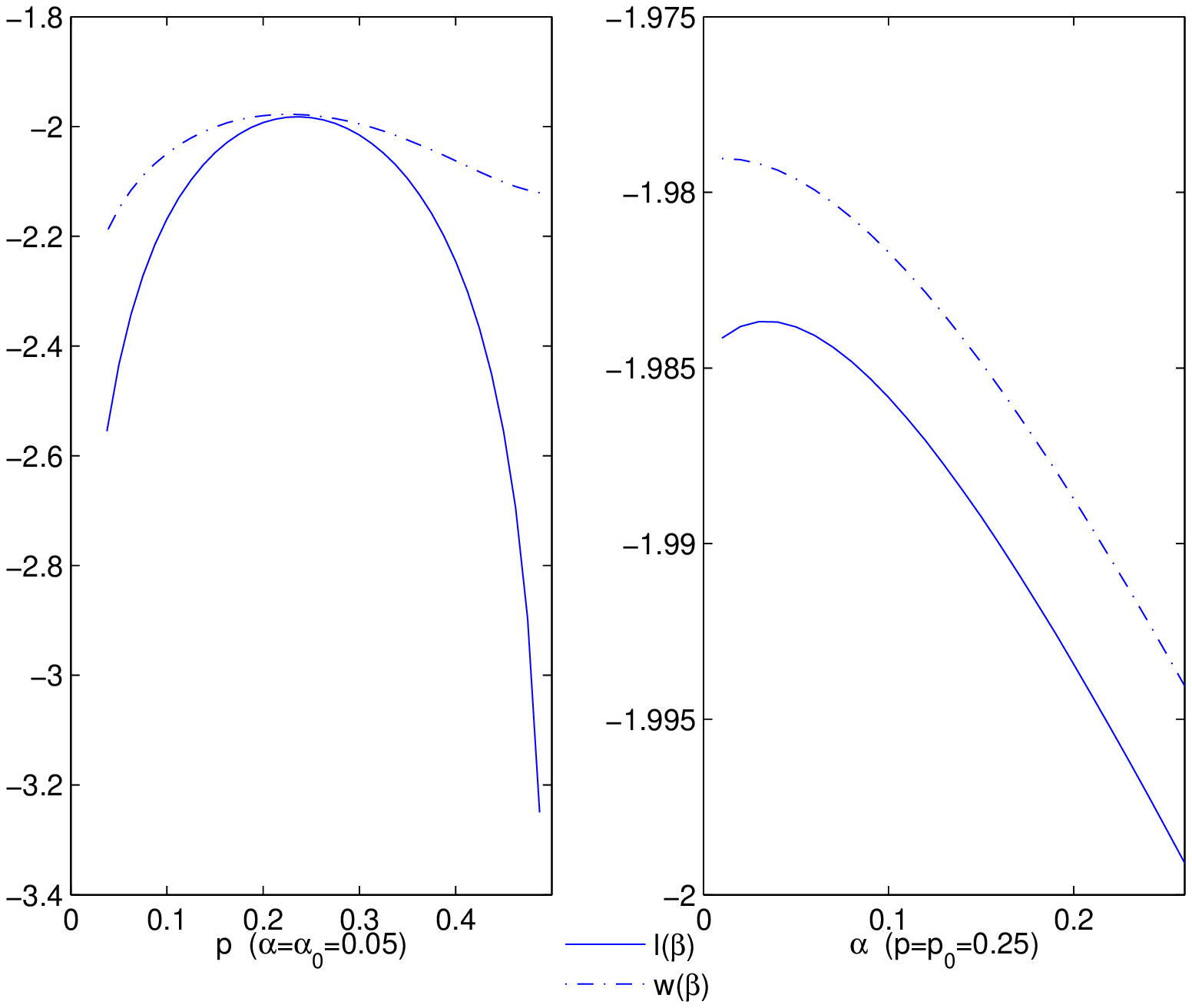}
 \caption{On top: $\ell$ and $w$ for the i.i.d. model
 ($p_0=0.25,\alpha_0=0.05$). On bottom: cuts of $\ell$ and $w$ for
 $\alpha=\alpha_0$ fixed and for $p=p_0$ fixed.}
\label{lqgraph}
\end{figure}

That can be explained by looking at the graph of $w(\beta)$ and
comparing it to $\ell(\beta)$ (Figure \ref{lqgraph}). We see that
both $w$ and $\ell$ are very flat with respect to $\alpha$ and as we deal with
numerical precision errors, finding out the true maximum value becomes
impossible. However, for $p=p_0$ if we look closely at the cuts of
$\ell$ and $w$ we appreciate that $\ell$ takes its maximum on
$\alpha_0$ and $w$ near this point. As the maximisation problem
complexity is reduced in this case we are able to find a quite good
estimation for $\alpha$. Concerning $p$, we see that both $\ell$ and $w$ have
a clear maximum near $p_0$, but again $\ell$ is less flat than $w$ at
this point. This is not surprising since  $\ell$ really is the
information divergence rate of the model.

\subsection{Simulations with Markov chains satisfying Assumption 2}
\begin{figure}[!h]
 \centering
\includegraphics[width=7.5cm]{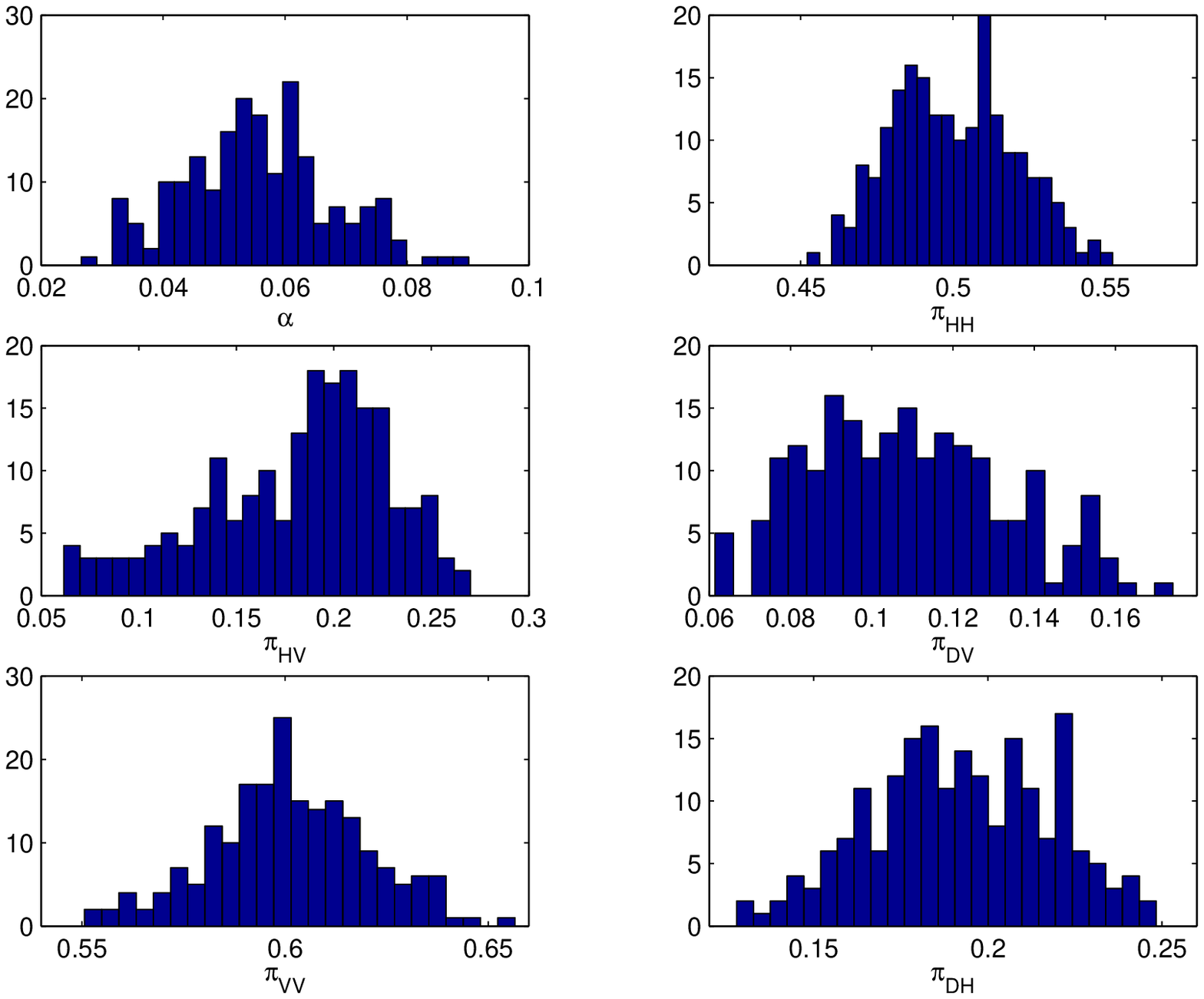}
 \includegraphics[width=7.5cm]{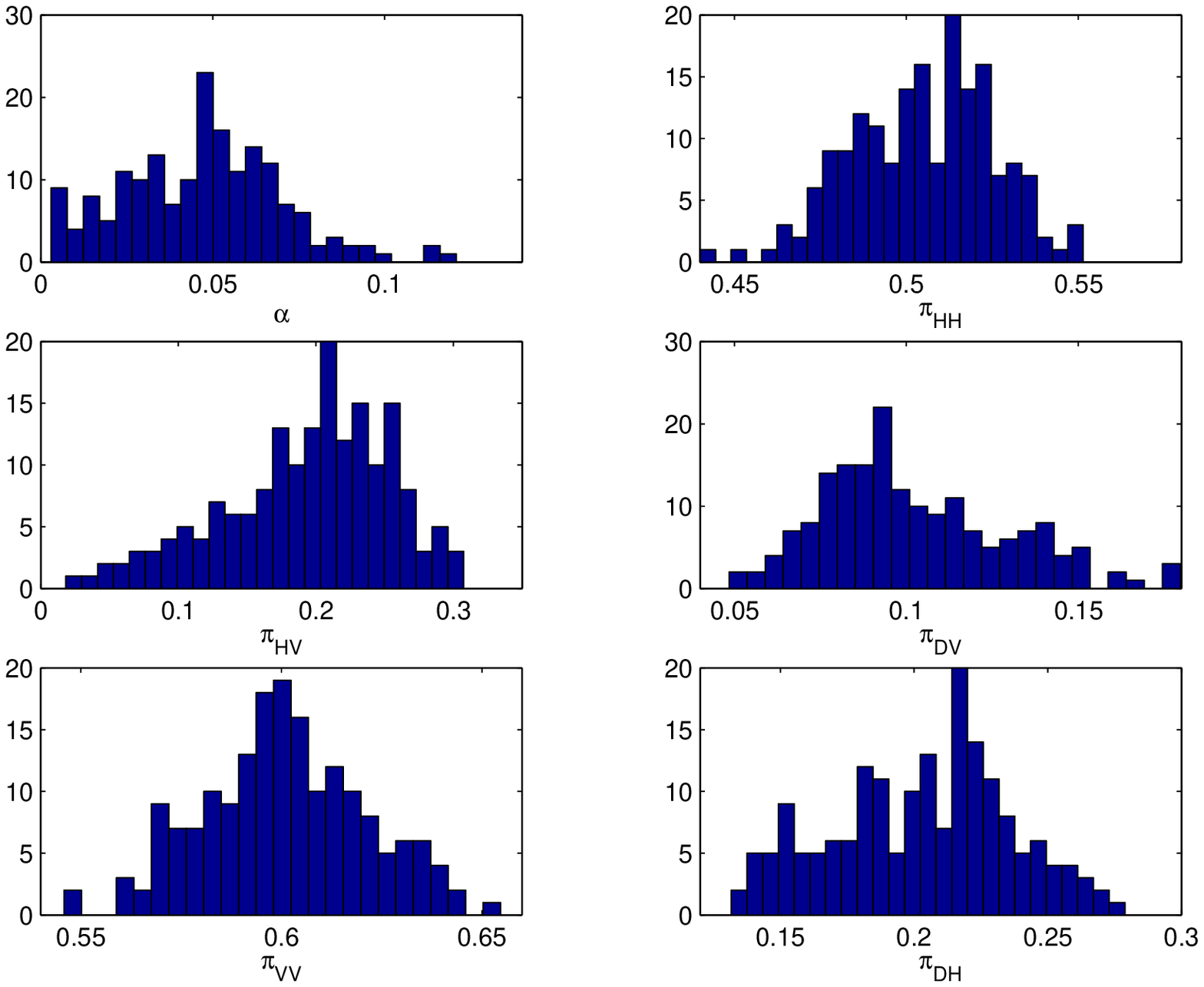}
 \caption{Histograms of maximum likelihood estimations of parameters obtained with  200 simulations
  from  the Markov chain model.
  On the left: estimation of the transition probabilities
  given $\alpha=\alpha_0$ and estimation of $\alpha$ given the true value of the transition
  probabilities. On the right: joint estimation of the transition probabilities and $\alpha$.}
 \label{markov}
\end{figure}

We have simulated 200 alignments of length 15000 with substitution
rate $\alpha_0=0.05$ and the following transition matrix for
$(\e_s)_s$\\
\begin{center}
$\begin{array}{ccc}\quad \quad D & \quad H & \quad V
\end{array}$\\
$\begin{array}{l}
D\\
H\\
V \end{array} \left(
\begin{array}{ccc}
0.7 & 0.2 & 0.1\\
0.3 & 0.5 & 0.2\\
0.3 & 0.1 & 0.6
\end{array}
\right)$
\end{center}
with initial distribution $p_0=q_0=0.25$. We have set as free
parameters $\pi_{HH}$, $\pi_{HV}$, $\pi_{DV}$, $\pi_{VV}$ and
$\pi_{DH}$. The equilibrium probability of every nucleotide is again
fixed to $0.25$. We can observe in Figure \ref{markov} that the maximum
likelihood estimators for these parameters and for $\alpha$ are
close to their true values even when the estimation is done
jointly. 
These results are rather encouraging since the Markov case is the
interesting one in biological applications.

% Catherine Matias, Laboratoire Statistique et G\'enome,  
%  Tour Evry  2, 523  pl.   des Terrasses  de l'Agora,  91~000 Evry,  France.
% E-mail:  matias@genopole.cnrs.fr

\bibliographystyle{chicago}
\bibliography{Phmm_final}

\begin{thebibliography}{}

\bibitem[\protect\citeauthoryear{Arribas~Gil, Metzler, and
  Plouhinec}{Arribas~Gil et~al.}{2005}]{Ana}
Arribas~Gil, A., D.~Metzler, and J.-L. Plouhinec (2005).
\newblock A fragment insertion and deletion model allowing fast and slow
  fragments.
\newblock Manuscript.

\bibitem[\protect\citeauthoryear{Baum and Petrie}{Baum and
  Petrie}{1966}]{BaumPetrie}
Baum, L. and T.~Petrie (1966).
\newblock Statistical inference for probabilistic functions of finite state
  {M}arkov chains.
\newblock {\em Ann. Math. Statist.\/} {\em 37}, pp.\  1554--1563.

\bibitem[\protect\citeauthoryear{Bishop and Thompson}{Bishop and
  Thompson}{1986}]{Bishop}
Bishop, M. and E.~Thompson (1986).
\newblock Maximum likelihood alignment of {DNA} sequences.
\newblock {\em J. Mol. Biol.\/} {\em 190}, pp.\  159--165.

\bibitem[\protect\citeauthoryear{Caliebe and Rösler}{Caliebe and
  Rösler}{2002}]{Caliebe}
Caliebe, A. and U.~Rösler (2002).
\newblock {Convergence of the maximum a posteriori path estimator in hidden
  Markov models.}
\newblock {\em IEEE Trans. Inf. Theory\/} {\em 48\/}(7), pp.\  1750--1758.

\bibitem[\protect\citeauthoryear{Cover and Thomas}{Cover and
  Thomas}{1991}]{Cover}
Cover, T.~M. and J.~A. Thomas (1991).
\newblock {\em {Elements of information theory.}}
\newblock {New York, USA}: {Wiley Series in Telecommunications, John Wiley \&
  Sons}.

\bibitem[\protect\citeauthoryear{Csiszár and Körner}{Csiszár and
  Körner}{1981}]{Csiszar}
Csiszár, I. and J.~Körner (1981).
\newblock {\em {Information theory. Coding theorems for discrete memoryless
  systems.}}
\newblock New York-San Francisco-London: {Probability and Mathematical
  Statistics. Academic Press.}

\bibitem[\protect\citeauthoryear{Davey and MacKay}{Davey and
  MacKay}{2001}]{Davey}
Davey, M.~C. and D.~J. MacKay (2001).
\newblock {Reliable communication over channels with insertions, deletions, and
  substitutions.}
\newblock {\em IEEE Trans. Inf. Theory\/} {\em 47\/}(2), pp.\  687--698.

\bibitem[\protect\citeauthoryear{Durbin, Eddy, Krogh, and Mitchison}{Durbin
  et~al.}{1998}]{Durbin}
Durbin, R., S.~Eddy, A.~Krogh, and G.~Mitchison (1998).
\newblock {\em {Biological sequence analysis: probabilistic models of proteins
  and nucleic acids}}.
\newblock {Cambridge, UK}: {Cambridge University Press}.

\bibitem[\protect\citeauthoryear{Hein, Wiuf, Knudsen, Moller, and Wibling}{Hein
  et~al.}{2000}]{HeinWiuf}
Hein, J., C.~Wiuf, B.~Knudsen, M.~Moller, and G.~Wibling (2000).
\newblock Statistical alignment: computational properties, homology testing and
  goodness-of-fit.
\newblock {\em J. Mol. Biol.\/} {\em 302}, pp.\  265--279.

\bibitem[\protect\citeauthoryear{Hobolth and Jensen}{Hobolth and
  Jensen}{2005}]{Hobolth}
Hobolth, A. and J.~Jensen (2005).
\newblock Applications of hidden {M}arkov models for characterization of
  homologous {DNA} sequences with a common gene.
\newblock {\em J. Comput. Biol.\/} {\em 12\/}(2), pp.\  186--203.

\bibitem[\protect\citeauthoryear{Holmes}{Holmes}{2005}]{Holmes}
Holmes, I. (2005).
\newblock Using evolutionary {E}xpectation {M}aximization to estimate indel
  rates.
\newblock {\em Bioinformatics\/} {\em 21\/}(10), pp.\  2294--2300.

\bibitem[\protect\citeauthoryear{Ibragimov and Has'minskii}{Ibragimov and
  Has'minskii}{1981}]{IH}
Ibragimov, I.~A. and R.~Z. Has'minskii (1981).
\newblock {\em Statistical Estimation. Asymptotic Theory}.
\newblock New York - Heidelberg -Berlin: Applications of Mathematics, Vol. 16.
  Springer-Verlag.

\bibitem[\protect\citeauthoryear{Kingman}{Kingman}{1968}]{Kingman}
Kingman, J. F.~C. (1968).
\newblock The ergodic theory of subadditive stochastic processes.
\newblock {\em J. Roy. Statist. Soc. Ser. B\/} {\em 30}, pp.\  499--510.

\bibitem[\protect\citeauthoryear{Knudsen and Miyamoto}{Knudsen and
  Miyamoto}{2003}]{KnudsenMiyamoto}
Knudsen, B. and M.~Miyamoto (2003).
\newblock Sequence alignments and pair hidden {M}arkov models using
  evolutionary history.
\newblock {\em J. Mol. Biol.\/} {\em 333}, pp.\  453--460.

\bibitem[\protect\citeauthoryear{Leroux}{Leroux}{1992}]{Leroux}
Leroux, B.~G. (1992).
\newblock Maximum-likelihood estimation for hidden {M}arkov models.
\newblock {\em Stochastic Process. Appl.\/} {\em 40\/}(1), pp.\  127--143.

\bibitem[\protect\citeauthoryear{Levenshtein}{Levenshtein}{2001}]{Levenshtein}
Levenshtein, V.~I. (2001).
\newblock {Efficient reconstruction of sequences.}
\newblock {\em IEEE Trans. Inf. Theory\/} {\em 47\/}(1), pp.\  2--22.

\bibitem[\protect\citeauthoryear{Metzler}{Metzler}{2003}]{Metzler}
Metzler, D. (2003).
\newblock Statistical alignment based on fragment insertion and deletion
  models.
\newblock {\em Bioinformatics\/} {\em 19\/}(4), pp.\  490--499.

\bibitem[\protect\citeauthoryear{Meyer and Durbin}{Meyer and
  Durbin}{2002}]{MeyerDurbin}
Meyer, I. and R.~Durbin (2002).
\newblock Comparative ab initio prediction of gene structures using pair
  {HMM}s.
\newblock {\em Bioinformatics\/} {\em 18\/}(10), pp.\  1309--1318.

\bibitem[\protect\citeauthoryear{Miklos, Lunter, and Holmes}{Miklos
  et~al.}{2004}]{Miklos}
Miklos, I., G.~A. Lunter, and I.~Holmes (2004).
\newblock {A "Long Indel" Model For Evolutionary Sequence Alignment}.
\newblock {\em Mol Biol Evol\/} {\em 21\/}(3), pp.\  529--540.

\bibitem[\protect\citeauthoryear{Needleman and Wunsch}{Needleman and
  Wunsch}{1970}]{Needleman}
Needleman, S. and C.~Wunsch (1970).
\newblock {A general method applicable to the search for similarities in the
  amino acid sequence of two proteins.}
\newblock {\em J. Mol. Biol.\/} {\em 48}, pp.\  443--453.

\bibitem[\protect\citeauthoryear{Pachter, Alexandersson, and Cawley}{Pachter
  et~al.}{2002}]{Pachter}
Pachter, L., M.~Alexandersson, and S.~Cawley (2002).
\newblock Applications of generalized pair hidden {M}arkov models to alignment
  and gene finding problems.
\newblock {\em J. Comput. Biol.\/} {\em 9\/}(2), pp.\  389--399.

\bibitem[\protect\citeauthoryear{Sucheston}{Sucheston}{1963}]{such}
Sucheston, L. (1963).
\newblock On mixing and the {Z}ero-{O}ne law.
\newblock {\em J. Math. Anal. and Appl.\/} {\em 6}, pp.\  447--456.

\bibitem[\protect\citeauthoryear{Thorne, Kishino, and Felsenstein}{Thorne
  et~al.}{1991}]{TKF}
Thorne, J., H.~Kishino, and J.~Felsenstein (1991).
\newblock An evolutionary model for maximum likelihood alignment of {DNA}
  sequences.
\newblock {\em J. Mol. Evol.\/} {\em 33}, pp.\  114--124.

\bibitem[\protect\citeauthoryear{Thorne, Kishino, and Felsenstein}{Thorne
  et~al.}{1992}]{TKF2}
Thorne, J., H.~Kishino, and J.~Felsenstein (1992).
\newblock Inching toward reality: an improved likelihood model of sequence
  evolution.
\newblock {\em J. Mol. Evol.\/} {\em 34}, pp.\  3--16.

\end{thebibliography}

 \end{document}